\documentclass[10pt,twoside]{amsart}
\usepackage[T1]{fontenc}
\usepackage{lmodern,microtype}
\usepackage{amssymb,amsmath,mathrsfs,mathtools,amsfonts,amsthm,amscd}
\usepackage[a4paper,top=1.5cm,bottom=1.5cm,left=2.5cm,right=2.5cm]{geometry}
\usepackage{xcolor}
\usepackage{cite}
\usepackage[hidelinks]{hyperref}
\usepackage{tikz-cd}

\newtheorem{Theorem}{Theorem}
\newtheorem{Definition}[Theorem]{Definition}

\newtheorem{Proposition}[Theorem]{Proposition}
\newtheorem{Lemma}[Theorem]{Lemma}

\newtheorem{Corollary}[Theorem]{Corollary}

\newtheorem{Example}[Theorem]{Example}

\newcommand{\Mor}{\mathrm{Mor}}
\newcommand{\Br}{\mathrm{Br}}

\newcommand{\Gpd}{\mathbf{Gpd}}
\newcommand{\Fun}{\mathrm{Fun}}
\newcommand{\LCS}{\mathrm{LCS}}
\newcommand{\Proj}{\mathrm{Proj}}
\newcommand{\id}{\mathrm{id}}

\DeclareMathOperator{\Spec}{Spec}

\DeclareMathOperator{\GL}{GL}

\DeclareMathOperator{\Aut}{Aut}

\newcommand{\PI}{\mathscr{PI}}
\newcommand{\AlgLaws}{\operatorname{Alg}}
\newcommand{\HH}{\operatorname{HH}}
\newcommand{\Der}{\operatorname{Der}}
\newcommand{\Perf}{\operatorname{Perf}}
\newcommand{\Az}{\operatorname{Az}}
\newcommand{\idG}{\operatorname{Id}^{G}}
\newcommand{\Tens}{\mathcal{T}}
\newcommand{\FreePI}{\mathscr F}
\newcommand{\cEnd}{\mathcal{E}nd}

\newcommand{\Lie}{\operatorname{Lie}}
\newcommand{\op}{\mathrm{op}}
\DeclareMathOperator{\rank}{rank}
\begin{document}

 \title[Moduli and K-Theory of Graded PI-Algebras]
{Moduli, Deformations and Algebraic K-Theory of Graded PI-Algebras}

\author{Lucio Centrone}
\address{Dipartimento di Matematica, Universit\`a degli Studi di Bari Aldo Moro, Via Edoardo Orabona, 4, 70125 Bari, Italy}
\address{IMECC, UNICAMP, Rua S\'ergio Buarque de Holanda 651, 13083-859 Campinas, SP, Brazil}
\email{lucio.centrone@uniba.it, centrone@unicamp.br}

\author{Maur\'icio Corr\^ea}
\address{Dipartimento di Matematica, Universit\`a degli Studi di Bari Aldo Moro, Via Edoardo Orabona, 4, 70125 Bari, Italy}
\email{mauricio.barros@uniba.it, mauricio.correa.mat@gmail.com}


\begin{abstract}
Over a field of characteristic zero, for a finite grading group and a finitely generated graded $T$-ideal, we construct relatively free sheaves and finite-rank quotient stacks of algebra laws with prescribed identities. Their tangent complexes are PI-restricted Hochschild complexes, and the second variations of the identities give explicit lifting obstructions. The first-order deformation groupoid carries a natural relative algebraic K-theory functor; for its square-zero extensions, the relative Chern character identifies rational relative K-theory with negative cyclic homology. Additional identities determine closed PI-strata and localisation fibre sequences in K-theory with Azumaya coefficients. On the Azumaya substack, the degree annihilates the Brauer class, and Morita transport yields a base-change-compatible K-theory of PI-Azumaya families together with a K-theoretic monodromy equivalence.
\end{abstract}

\maketitle

\section{Introduction}\label{sec:introduction}

Noncommutative geometry treats a noncommutative algebra as an algebra of functions on a space not directly accessible through ordinary points.  This principle underlies Connes' differential-geometric approach, formal noncommutative geometry and representation-theoretic methods in PI-theory \cite{Connes,Kontsevich,Kontsevich-Rosenberg}.  Graded geometry supplies a second source of examples: supergeometry encodes commutation by a grading and a sign bicharacter \cite{BerLe,Kostant,Manin,Deligne,DoWi,BHPo}, and $\mathbb Z_2^n$-supergeometry extends the same mechanism to higher colour gradings \cite{CGN,vysoky}.

Polynomial identities provide a common algebraic language for commutative and supercommutative algebras, finite-dimensional and matrix algebras, upper triangular and Clifford algebras, Azumaya algebras, algebras finite over their centres, quantum tori and several quantum algebras at roots of unity.  Artin's use of identities and invariant theory for Azumaya algebras already exhibited their geometric content \cite{Artin}; the modern theory relates identities, central localisation, generic matrices and representation spaces \cite{giambruno2005polynomial,drensky2000free,AGPR,ProcesiPI,RowenPI}.

Stacks of algebra structures and Hochschild-theoretic deformation theory have substantial antecedents \cite{LieblichCoherent,BehrendNoohi,HinichHomotopy,KontsevichSoibelmanOperads,YalinModuli}.  Our results concern instead the finite-rank geometry determined by a fixed finitely generated graded $T$-ideal: its closed law stack, the normal morphism obtained by linearising additional identities, and a second-order lifting criterion combining the Gerstenhaber quadratic term with the Hessians of those identities.  We relate these constructions to relatively free sheaves, Azumaya forms, algebraic K-theory and Morita-variable monodromy.

The breadth of these examples makes the PI language a natural point of departure, but one distinction is essential.  The equational class
\[
\operatorname{Var}^G(I)=\{B: I\subseteq\operatorname{Id}^G(B)\}
\]
is a category of algebras, not itself a geometric space.  An inclusion of identity ideals does not produce homomorphisms between arbitrary algebras; a matrix algebra over a local ring is usually not a local ring; and the derivations of an associative algebra form a Lie algebra rather than an associative algebra satisfying the same polynomial identities.  The geometry must consequently be built from the identities by constructions that retain families, changes of basis, descent and infinitesimal variation.

We develop this programme at two complementary levels.

\smallskip
\noindent\emph{The sheaf-theoretic level.}
On a ringed space $(X,\mathcal O_X)$ we construct the relatively free sheaf
\[
\FreePI^G_{I,X}(\mathcal E)=\Tens_X^G(\mathcal E)/\mathcal J_I(\mathcal E).
\]
Its universal, stalk and pullback properties are proved in Theorem~\ref{thm:sheaf-relative-free} and Proposition~\ref{prop:free-restriction-diagram}.  The centrally locally ringed PI-spaces of Definition~\ref{def:central-PI-space} retain the commutative local geometry in the centre while allowing matrix, Azumaya and Clifford coefficient sheaves whose noncommutative stalks need not be local rings.

\smallskip
\noindent\emph{The moduli level.}
For a finite graded dimension vector $\mathbf d$, Definition~\ref{def:coordinate-ring-law-scheme} gives the affine law scheme $\AlgLaws^G_{\mathbf d}(I)$.  The quotient
\[
\PI^G_{I,\mathbf d}=[\AlgLaws^G_{\mathbf d}(I)/H_{\mathbf d}]
\]
is proved algebraic in Theorem~\ref{thm:PI-stack-algebraicity} and is identified with the groupoid of finite locally free graded algebra families satisfying $I$ in Theorem~\ref{thm:PI-stack-classification}.  Thus the identities are the equations of the moduli problem, not labels attached after the fact.

The two levels meet through the classifying map of Proposition~\ref{prop:sheaf-stack-interface} and the tautological family of Proposition~\ref{prop:tautological-family}.  For $I\subseteq J$, additional identities induce both the quotient
\[
\FreePI^G_{I,X}(\mathcal E)\twoheadrightarrow\FreePI^G_{J,X}(\mathcal E)
\]
of Corollary~\ref{cor:sheaf-identity-quotient} and the closed immersion
\[
\PI^G_{J,\mathbf d}\hookrightarrow\PI^G_{I,\mathbf d}
\]
of Theorem~\ref{thm:closed-PI-immersion}.
The infinitesimal geometry is PI-specific.  Theorem~\ref{thm:tangent-complex} identifies the non-positive truncation of the intrinsic tangent complex at $[B]$ with
\[
C^1_{G,\mathrm n}(B,B)\xrightarrow{d_{\mathrm{Hoch}}}Z^2_{G,I,\mathrm n}(B,B),
\]
and identifies its associated groupoid with first-order deformations preserving $I$.  For $I\subseteq J$, Theorem~\ref{thm:relative-PI-tangent} identifies the first-order normal quotient with $Z^2_{G,I,\mathrm n}(B,B)/Z^2_{G,J,\mathrm n}(B,B)$ and computes the normal map by linearising the additional identities.  Proposition~\ref{prop:associative-harrison} recovers Hochschild and Harrison cohomology in the associative and commutative cases.

First-order compatibility does not imply liftability.  Theorem~\ref{thm:second-order-PI-obstruction} computes the obstruction to extending $\mu+t\varphi$ over $k[t]/(t^3)$ from the Gerstenhaber quadratic term and the second variations of generators of $I$; its vanishing is equivalent to the existence of the lift.  Corollary~\ref{cor:smooth-formal-liftability} gives formal liftability at smooth points of the law scheme.

The deformation theory has a K-theoretic realisation.  Theorem~\ref{thm:PI-relative-K-functor} assigns a relative K-theory spectrum functorially to every PI-preserving first-order deformation, and Theorem~\ref{thm:Goodwillie-PI-relative} identifies its positive rational homotopy groups with relative negative cyclic homology.  Example~\ref{ex:dual-numbers-relative-K} computes the first relative K-group for a nontrivial deformation of the dual numbers.  For $I\subseteq J$, Theorem~\ref{thm:PI-identity-localisation} gives a localisation fibre sequence for the nonconnective K-theory of the closed Azumaya stratum defined by the additional identities.

The Azumaya locus is open and stable under base change by Theorem~\ref{thm:Azumaya-open}; its derivations and form substacks are described in Theorem~\ref{thm:Azumaya-derivations} and Proposition~\ref{prop:orbit-classifying-stack}.  If $N=\sum_gd_g$, Theorem~\ref{thm:brauer-torsion-on-PI-stack} shows that the locus is empty unless $N=n^2$ and that every rank-$n^2$ family has Brauer class annihilated by $n$.  Theorem~\ref{thm:PI-Azumaya-K-pseudofunctor} organises its K-theory by the pseudofunctor
\[
\mathbf K_{I,\mathbf d,S}:\bigl(\PI^{G,\Az}_{I,\mathbf d}(S)\bigr)^{\Mor}\longrightarrow\mathbf{Sp}^{\simeq},
\qquad \mathcal A\longmapsto\mathbf K(S;\mathcal A),
\]
obtained by restricting Morita-invariant K-theory from degree-$n$ Azumaya algebras.  Its equivalence type depends only on the Brauer class, and it is compatible with pullback in $S$.  Example~\ref{ex:matrix-forms-K} identifies the matrix identity locus with $B\operatorname{PGL}_n$ and sends a torsor to the K-theory spectrum of its associated Azumaya algebra.  Corollary~\ref{cor:K-Azumaya-after-degree} compares this twisted K-theory with that of the base after $n$ is inverted.

The final part concerns variable coefficients.  Proposition~\ref{prop:Morita-module-transport} constructs transport along invertible bimodules.  The fixed-coefficient equivalence
\[
\LCS(X,B)\simeq\Fun(\Pi_1(X),\Proj(B)^{\simeq})
\]
is Theorem~\ref{thm:fixed-coefficient-monodromy}; Theorem~\ref{thm:expanded-Morita-monodromy} proves its pseudonaturality over a specified Morita $(2,1)$-groupoid and the resulting biequivalence of Grothendieck totals.  The exact enhancement of Theorem~\ref{thm:K-monodromy} induces equivalences of algebraic K-theory spectra, coherently compatible with Morita bimodules.  These results give a Morita-equivariant Betti classification without invoking flat connections or $\mathcal D$-modules.

The paper is organised as follows.  Section~\ref{sec:graded-relfree} develops graded identities and relatively free algebras; Section~\ref{sec:pi-sheaves} treats their sheaf-theoretic form.  We construct the finite-rank moduli problem in Section~\ref{sec:pi-moduli}.  Section~\ref{sec:pi-deformations} develops PI-restricted deformation theory, relative K-theory, the relative Chern character and localisation along additional identities.  Section~\ref{sec:azumaya} studies Azumaya loci, Brauer torsion and their K-theory.  Morita transport and variable-coefficient monodromy, including the exact K-theoretic enhancement, occupy Sections~\ref{sec:morita-transport} and~\ref{sec:morita-monodromy}.  Section~\ref{sec:examples} collects the examples.

\subsection*{Acknowledgements}
MC is partially supported by the Universit\`a degli Studi di Bari and by the PRIN 2022MWPMAB ``Interactions between Geometric Structures and Function Theories'', and is a member of INdAM-GNSAGA; he is also partially supported by Fapemig grants APQ-02674-21, APQ-00798-18, and APQ-00056-20. LC is partially supported by the Universit\`a degli Studi di Bari, INdAM-GNSAGA, and the PNRR MUR project PE0000023-NQSTI.
\section{Graded identities and relatively free algebras}\label{sec:graded-relfree}

Throughout the paper, $F$ is a field of characteristic zero and $G$ is a finite group with neutral element $e$.  All algebras are associative and unital, and all homomorphisms preserve the unit.

\begin{Definition}
A $G$-grading on an $F$-algebra $B$ is a decomposition $B=\bigoplus_{g\in G}B_g$ such that $B_gB_h\subseteq B_{gh}$.  A homomorphism is graded if it preserves every homogeneous component.
\end{Definition}

Let $F\langle X\mid G\rangle$ be the free associative algebra on variables $x_i^g$ of prescribed degrees $g\in G$.  A graded polynomial $f$ is a graded identity of $B$ if every degree-preserving evaluation of $f$ in $B$ is zero.  The set $\idG(B)$ is a graded $T$-ideal, namely an ideal stable under every graded endomorphism of the free algebra.

\begin{Definition}\label{def:graded-equational-category}
For a graded $T$-ideal $I$, let $G\text{-}\mathbf{Alg}_I$ be the category whose objects are unital $G$-graded $F$-algebras $B$ with $I\subseteq\idG(B)$ and whose morphisms are unital graded homomorphisms.  We write
\[
\operatorname{Var}^G(I)=\operatorname{Ob}(G\text{-}\mathbf{Alg}_I).
\]
If $I=\idG(A)$, we also write $\operatorname{Var}^G(A)$.  For a graded vector space $V$, let $T^G(V)$ be its graded tensor algebra and let
\[
J_I(V)=\sum_{\rho}\,T^G(V)\,\rho(I)\,T^G(V),
\]
where $\rho$ runs over all graded algebra homomorphisms from $F\langle X\mid G\rangle$ to $T^G(V)$.  The relatively free algebra of the theory $I$ on $V$ is
$
\FreePI_I^G(V):=T^G(V)/J_I(V).
$
\end{Definition}

The definition is independent of a homogeneous basis of $V$ and has the expected universal property.

\begin{Proposition}\label{prop:relative-free-adjunction-vector}
For every graded algebra $B$ satisfying $I$, restriction to $V$ induces a natural bijection
\[
\operatorname{Hom}_{G\text{-}\mathrm{Alg},I}(\FreePI_I^G(V),B)
\cong
\operatorname{Hom}_{G\text{-}\mathrm{Vect}}(V,B).
\]
Consequently, $\FreePI_I^G$ is left adjoint to the forgetful functor from graded algebras satisfying $I$ to graded vector spaces.
\end{Proposition}

\begin{proof}
A graded linear map $V\to B$ extends uniquely to a graded algebra map $T^G(V)\to B$.  Since $B$ satisfies $I$, every evaluation generating $J_I(V)$ is sent to zero, so the map factors uniquely through $\FreePI_I^G(V)$.  The converse is restriction to the canonical image of $V$.
\end{proof}

\subsection{The adjunction, unit and counit}
Let
$
U_I:G\text{-}\mathbf{Alg}_I\longrightarrow G\text{-}\mathbf{Vect}_F
$
be the forgetful functor. Proposition~\ref{prop:relative-free-adjunction-vector} is represented by the adjunction
$
\FreePI_I^G\dashv U_I.
$
Its unit at a graded vector space $V$ is the canonical inclusion
\[
\eta_V:V\longrightarrow U_I\FreePI_I^G(V),
\]
and its counit at $B\in G\text{-}\mathbf{Alg}_I$ is the evaluation morphism
$
\epsilon_B:\FreePI_I^G(U_I B)\longrightarrow B.
$ 
For every graded linear map $u:V\to U_I B$, the corresponding algebra morphism $\bar u$ is characterised by
\[
\begin{tikzcd}[column sep=large,row sep=large]
V \arrow[r,"\eta_V"] \arrow[dr,"u"']
& U_I\FreePI_I^G(V) \arrow[d,"U_I(\bar u)"]\\
& U_I B.
\end{tikzcd}
\]
The triangular identities are
$
\epsilon_{\FreePI_I^G(V)}\circ\FreePI_I^G(\eta_V)=\id_{\FreePI_I^G(V)}
 $ and $
U_I(\epsilon_B)\circ\eta_{U_I B}=\id_{U_I B}.
$
The adjunction adjoins multiplication subject precisely to the identities in $I$.

\begin{Proposition}\label{prop:relative-free-pushout}
Let $q_I:T^G(V)\twoheadrightarrow\FreePI_I^G(V)$ be the quotient. For every $B\in G\text{-}\mathbf{Alg}_I$ and graded linear map $u:V\to B$, the square
\[
\begin{tikzcd}[column sep=large,row sep=large]
T^G(V) \arrow[r,"q_I"] \arrow[d,"T^G(u)"']
& \FreePI_I^G(V) \arrow[d,"\bar u"]\\
T^G(U_I B) \arrow[r,"\operatorname{ev}_B"']
& B
\end{tikzcd}
\]
commutes, and $\bar u$ is the unique lower-right arrow with this property. Equivalently, $q_I$ is initial among graded algebra quotients of $T^G(V)$ whose targets satisfy $I$.
\end{Proposition}
\begin{proof}
Here $\operatorname{ev}_B:T^G(U_I B)\to B$ is tensor-word evaluation; it factors as $T^G(U_I B)\twoheadrightarrow\FreePI_I^G(U_I B)\xrightarrow{\epsilon_B}B$.  The left-hand composite evaluates every tensor word in the chosen elements $u(V)$. Since $B$ satisfies $I$, it kills $J_I(V)$ and hence factors uniquely through $q_I$. The displayed square expresses this factorisation.
\end{proof}

Finite generation of graded identity ideals is essential in what follows.

\begin{Theorem}\label{thm:graded-specht}
Let $A$ be an associative PI-algebra endowed with a $G$-grading.  Then $\idG(A)$ is finitely generated as a graded $T$-ideal.
\end{Theorem}

\begin{proof}
This is the graded Specht theorem of Aljadeff and Kanel-Belov for finite $G$ in characteristic zero; see \cite{AljadeffBelov2010}.
\end{proof}

By polarisation, a finite generating system may be replaced by finitely many multihomogeneous multilinear generators.  We fix such a system whenever a finite-type presentation is needed.

\begin{Example}
A commutative algebra satisfies $[x_1,x_2]=x_1x_2-x_2x_1$.  Every finite-dimensional algebra is PI, because every alternating multilinear polynomial in more variables than its dimension vanishes.  The infinite Grassmann algebra is PI and, with its canonical $\mathbb Z_2$-grading, satisfies the graded supercommutativity identities.
\end{Example}

\subsection{Equational closure and the role of the \texorpdfstring{$T_G$}{T-G}-ideal}
The notation $\operatorname{Var}^G(I)$ is justified by Birkhoff's equational viewpoint.  The class of graded algebras satisfying $I$ is closed under graded subalgebras, graded quotients and arbitrary products.  Conversely, a class closed under these operations and under graded isomorphism is controlled by the graded identities common to its objects.  Only the forward implication is needed here; it also explains why $I$, rather than a chosen presentation of an algebra, is the natural parameter.

\begin{Proposition}\label{prop:equational-closure}
Let $I$ be a graded $T$-ideal.
\begin{enumerate}
\item Graded subalgebras, graded quotients and products of algebras in $\operatorname{Var}^G(I)$ again belong to $\operatorname{Var}^G(I)$.
\item Filtered colimits of graded algebras in $\operatorname{Var}^G(I)$ belong to $\operatorname{Var}^G(I)$.
\item If $I\subseteq J$, then $\operatorname{Var}^G(J)\subseteq\operatorname{Var}^G(I)$.
\end{enumerate}
\end{Proposition}
\begin{proof}
Polynomial evaluation is compatible with subalgebras, quotient maps and coordinatewise products.  A polynomial involves only finitely many elements, so an evaluation in a filtered colimit is represented at a finite stage.  The final assertion follows directly from the definition.
\end{proof}

The reversal reflects the geometric fact that additional identities cut out a closed substack.

\subsection{Basic models}
\begin{Example}
For the trivial grading, commutativity is the identity $[x_1,x_2]=x_1x_2-x_2x_1$.  For a bicharacter $\varepsilon:G\times G\to F^*$, the identities $x^gy^h-\varepsilon(g,h)y^hx^g$ define the $\varepsilon$-commutative, or colour-commutative, variety.  Supercommutativity is the case $G=\mathbb Z_2$ and $\varepsilon(i,j)=(-1)^{ij}$.
\end{Example}

\begin{Example}
If $B$ has dimension $d$, every alternating multilinear polynomial in $d+1$ variables vanishes on $B$.  Thus matrix, triangular, Clifford and finite-dimensional Hopf algebras are PI\@.  The identity ideal usually contains much more information than the alternating identity; for instance, the Amitsur--Levitzki identity detects the matrix degree.
\end{Example}

\begin{Example}
Let $E$ be the infinite Grassmann algebra on generators $e_i$ with $e_ie_j=-e_je_i$.  As an ungraded algebra it satisfies the triple commutator identity $[x_1,x_2,x_3]=0$.  With the canonical $\mathbb Z_2$-grading, it is the relatively free model for supercommutative identities in characteristic zero.  More refined $\mathbb Z$-graded identities are described in \cite{alanplamen}.
\end{Example}

\subsection{Functoriality of the relatively free algebra}
\begin{Proposition}\label{prop:relative-free-two-variable-functoriality}
The assignments
\[
V\longmapsto\FreePI_I^G(V),\qquad
I\longmapsto\FreePI_I^G(V)
\]
have the following functoriality.
\begin{enumerate}
\item A graded linear map $u:V\to W$ induces a unique graded algebra homomorphism
\[
\FreePI_I^G(u):\FreePI_I^G(V)\longrightarrow\FreePI_I^G(W)
\]
whose restriction to $V$ is the composite $V\xrightarrow{u}W\to\FreePI_I^G(W)$.
\item If $I\subseteq J$, then the identity map of $T^G(V)$ induces a natural surjection
$
q_{I,J,V}:\FreePI_I^G(V)\twoheadrightarrow\FreePI_J^G(V)
$
with
\[
\ker q_{I,J,V}=J_J(V)/J_I(V).
\]
\item For $I\subseteq J\subseteq K$ and $u:V\to W$ one has
$
q_{J,K,W}\circ q_{I,J,W}\circ\FreePI_I^G(u)
=
\FreePI_K^G(u)\circ q_{I,K,V}.
$
\end{enumerate}
Thus $(I,V)\mapsto\FreePI_I^G(V)$ is contravariant in the poset of identity ideals and covariant in graded vector spaces.
\end{Proposition}
\begin{proof}
The first assertion follows from the adjunction in Proposition~\ref{prop:relative-free-adjunction-vector}.  If $I\subseteq J$, then every evaluation of an element of $I$ is an evaluation of an element of $J$, so $J_I(V)\subseteq J_J(V)$; the second assertion is the third isomorphism theorem.  The final identity follows because all maps are induced by the same tensor-algebra map and the inclusions of the corresponding relation ideals.
\end{proof}
\section{Relatively free PI-algebra sheaves}\label{sec:pi-sheaves}

Let $(X,\mathcal O_X)$ be a ringed space over $F$.  A graded $\mathcal O_X$-module is a direct sum $\mathcal E=\bigoplus_{g\in G}\mathcal E_g$.  Tensor products below are tensor products of sheaves; in particular, they are sheafifications of the corresponding presheaf tensor products.

\begin{Definition}
The graded tensor algebra of $\mathcal E$ is
\[
\Tens_X^G(\mathcal E)=\bigoplus_{m\geq0}\mathcal E^{\otimes_{\mathcal O_X}m},
\]
with the grading induced by multiplication in $G$.  For a graded $T$-ideal $I$, let $\mathcal J_I(\mathcal E)$ be the smallest sheaf of two-sided ideals of $\Tens_X^G(\mathcal E)$ which contains
$
\rho(f) $ with $ f\in I,
$
for every local graded algebra homomorphism
\[
\rho:F\langle X\mid G\rangle\longrightarrow\Tens_X^G(\mathcal E)|_U.
\]
Equivalently, it is the sheafification of the presheaf of ideals generated by all degree-preserving evaluations of elements of $I$ on homogeneous local sections of the tensor algebra.  Set
\[
\FreePI_{I,X}^G(\mathcal E):=
\Tens_X^G(\mathcal E)/\mathcal J_I(\mathcal E).
\]
\end{Definition}

\begin{Theorem}\label{thm:sheaf-relative-free}
Let $\mathcal E$ be a graded $\mathcal O_X$-module.
\begin{enumerate}
\item The sheaf $\FreePI_{I,X}^G(\mathcal E)$ is a sheaf of unital graded $\mathcal O_X$-algebras satisfying $I$.
\item For every sheaf $\mathcal B$ of graded $\mathcal O_X$-algebras satisfying $I$, restriction to $\mathcal E$ gives a natural bijection
\[
\operatorname{Hom}_{G\text{-}\mathcal O_X\mathrm{Alg},I}
(\FreePI_{I,X}^G(\mathcal E),\mathcal B)
\cong
\operatorname{Hom}_{G\text{-}\mathcal O_X\mathrm{Mod}}(\mathcal E,\mathcal B).
\]
\item For every $x\in X$ there is a canonical isomorphism
$
(\FreePI_{I,X}^G(\mathcal E))_x
\cong
\FreePI_{I,\mathcal O_{X,x}}^G(\mathcal E_x).
$
\item If $f:Y\to X$ is a morphism of schemes, $\mathcal E$ is quasi-coherent, and $I$ is generated by finitely many multihomogeneous multilinear identities, then the canonical base-change morphism is an isomorphism
\[
f^*\FreePI_{I,X}^G(\mathcal E)
\xrightarrow{\sim}
\FreePI_{I,Y}^G(f^*\mathcal E).
\]
\end{enumerate}
\end{Theorem}

\begin{proof}
Put
$
\mathcal T=\Tens_X^G(\mathcal E), 
\mathcal J=\mathcal J_I(\mathcal E), 
\mathcal F=\mathcal T/\mathcal J,
$
and denote the quotient map by $q:\mathcal T\twoheadrightarrow\mathcal F$.
The quotient is taken in sheaves of graded $\mathcal O_X$-algebras, so
associativity, the grading, and the unit descend from $\mathcal T$.
We first prove that $\mathcal F$ satisfies $I$.  Let $U\subseteq X$ be
open, let $f\in I$, and let $\bar a_1,\ldots,\bar a_m$ be homogeneous
sections of $\mathcal F(U)$ of the degrees prescribed by $f$.  An
epimorphism of sheaves is locally surjective; hence there is an open cover
$U=\bigcup_\alpha U_\alpha$ and homogeneous lifts
$a_{i,\alpha}\in\mathcal T(U_\alpha)$ of
$\bar a_i|_{U_\alpha}$.  By the definition of $\mathcal J$,
\[
f(a_{1,\alpha},\ldots,a_{m,\alpha})\in\mathcal J(U_\alpha).
\]
Therefore
$
f(\bar a_1,\ldots,\bar a_m)|_{U_\alpha}
=q\bigl(f(a_{1,\alpha},\ldots,a_{m,\alpha})\bigr)=0
$
for every $\alpha$.  The sheaf axiom then gives
$f(\bar a_1,\ldots,\bar a_m)=0$ on $U$.  This proves (1).

For (2), let $u:\mathcal E\to\mathcal B$ be a graded
$\mathcal O_X$-linear map.  The universal property of the tensor algebra in
sheaves gives a unique graded algebra morphism
$
\widetilde u:\mathcal T\longrightarrow\mathcal B
$
whose restriction to $\mathcal E$ is $u$.  For every local graded
substitution $\rho$ and every $f\in I$, the section
$\widetilde u(\rho(f))$ is the corresponding evaluation of $f$ in
$\mathcal B$, hence is zero.  Thus $\widetilde u$ kills $\mathcal J$ and
factors uniquely through a morphism $\bar u:\mathcal F\to\mathcal B$.
Conversely, restriction of an algebra morphism $\mathcal F\to\mathcal B$
to $\mathcal E$ recovers $u$.  The two operations are inverse and are
functorial in both variables.

For (3), exactness of the stalk functor gives
$
\mathcal F_x\cong\mathcal T_x/\mathcal J_x.
$
Stalks commute with direct sums and tensor products, so there is a canonical
isomorphism
\[
\mathcal T_x\cong T^G_{\mathcal O_{X,x}}(\mathcal E_x).
\]
It remains to identify the relation ideal.  A germ of a local evaluation of
an element of $I$ is an evaluation on germs, and therefore
$\mathcal J_x\subseteq J_I(\mathcal E_x)$.  Conversely, an evaluation in
$T^G_{\mathcal O_{X,x}}(\mathcal E_x)$ involves only finitely many germs.
Choose representatives on a common neighbourhood of $x$; the resulting
local evaluation belongs to $\mathcal J$, and its germ is the prescribed
element.  Hence $\mathcal J_x=J_I(\mathcal E_x)$, proving the stalk formula.

For (4), pullback of modules is right exact, commutes with arbitrary direct
sums, and has the canonical monoidal isomorphisms
\[
f^*(\mathcal M\otimes_{\mathcal O_X}\mathcal N)
\cong f^*\mathcal M\otimes_{\mathcal O_Y}f^*\mathcal N.
\]
Consequently
$
f^*\mathcal T\cong\Tens_Y^G(f^*\mathcal E).
$
Let $\mathcal K$ be the image of $f^*\mathcal J$ under this isomorphism.
Pullback carries an evaluation of a generator $f_j$ to the corresponding
evaluation after pullback, so
$\mathcal K\subseteq\mathcal J_I(f^*\mathcal E)$.  For the reverse
inclusion, work locally on $Y$.  Every section of a pullback module is a
finite sum of sections of the form $a\,f^{-1}s$.  Since each $f_j$ is
multilinear, an evaluation of $f_j$ on such sums expands as an
$\mathcal O_Y$-linear combination of pullbacks of evaluations of $f_j$ on
local sections of $\mathcal T$.  It therefore lies in $\mathcal K$.
Because the $f_j$ generate $I$ as a graded $T$-ideal, this proves
$\mathcal K=\mathcal J_I(f^*\mathcal E)$.  Applying $f^*$ to
$\mathcal T\twoheadrightarrow\mathcal F$ and using right exactness now
gives the asserted base-change isomorphism.
\end{proof}

\begin{Corollary}\label{cor:sheaf-identity-quotient}
If $I\subseteq J$ are graded $T$-ideals, then there is a natural surjection
$
\FreePI_{I,X}^G(\mathcal E)
\twoheadrightarrow
\FreePI_{J,X}^G(\mathcal E)
$
with kernel $\mathcal J_J(\mathcal E)/\mathcal J_I(\mathcal E)$.
\end{Corollary}

\begin{proof}
The inclusion $I\subseteq J$ gives $\mathcal J_I(\mathcal E)\subseteq\mathcal J_J(\mathcal E)$, and the assertion follows from the third isomorphism theorem.
\end{proof}

\begin{Proposition}\label{prop:tensor-constant-finite}
If $A$ is a finite-dimensional graded $F$-algebra, then
$
\mathcal A:=A\otimes_F\mathcal O_X
$
is a locally free graded $\mathcal O_X$-algebra and satisfies $\idG(A)$.  Its stalk at $x$ is $A\otimes_F\mathcal O_{X,x}$.
\end{Proposition}

\begin{proof}
Choose a homogeneous basis of $A$.  As an $\mathcal O_X$-module, $\mathcal A$ is a finite direct sum of copies of $\mathcal O_X$; the sheaf condition and the stalk formula therefore follow from this decomposition.  Every graded identity of $A$ remains an identity after extension of scalars by a commutative $F$-algebra.
\end{proof}

\subsection{Exact diagrams for stalks and pullback}
The construction of Theorem~\ref{thm:sheaf-relative-free} is controlled by
\begin{equation}\label{eq:sheaf-envelope-exact}
0\longrightarrow\mathcal J_I(\mathcal E)
\longrightarrow\Tens_X^G(\mathcal E)
\longrightarrow\FreePI_{I,X}^G(\mathcal E)
\longrightarrow0.
\end{equation}
Taking the stalk at $x$ gives a diagram with exact rows
\[
\begin{tikzcd}[column sep=small]
0 \arrow[r]
& \mathcal J_I(\mathcal E)_x \arrow[r] \arrow[d,"\sim"']
& \Tens_X^G(\mathcal E)_x \arrow[r] \arrow[d,"\sim"']
& \FreePI_{I,X}^G(\mathcal E)_x \arrow[r] \arrow[d,dashed,"\sim"]
& 0\\
0 \arrow[r]
& J_I(\mathcal E_x) \arrow[r]
& T^G_{\mathcal O_{X,x}}(\mathcal E_x) \arrow[r]
& \FreePI^G_{I,\mathcal O_{X,x}}(\mathcal E_x) \arrow[r]
& 0.
\end{tikzcd}
\]
The dashed isomorphism is induced on cokernels.  For a morphism $f:Y\to X$, the base-change morphism is the unique arrow completing
\[
\begin{tikzcd}[column sep=large,row sep=large]
f^*\Tens_X^G(\mathcal E) \arrow[r,two heads] \arrow[d,"\sim"']
& f^*\FreePI^G_{I,X}(\mathcal E) \arrow[d,dashed]\\
\Tens_Y^G(f^*\mathcal E) \arrow[r,two heads]
& \FreePI^G_{I,Y}(f^*\mathcal E).
\end{tikzcd}
\]
The finite multihomogeneous generating system for $I$ identifies the pullback of the upper relation sheaf with the lower relation sheaf, so the dashed arrow is an isomorphism.

\begin{Proposition}\label{prop:free-restriction-diagram}
For every open immersion $j:U\hookrightarrow X$ there is a canonical isomorphism
\[
j^{-1}\FreePI^G_{I,X}(\mathcal E)
\xrightarrow{\sim}
\FreePI^G_{I,U}(j^{-1}\mathcal E),
\]
characterised by
\[
\begin{tikzcd}[column sep=large,row sep=large]
j^{-1}\mathcal E \arrow[r] \arrow[d,equal]
& j^{-1}\FreePI^G_{I,X}(\mathcal E) \arrow[d,"\sim"]\\
j^{-1}\mathcal E \arrow[r]
& \FreePI^G_{I,U}(j^{-1}\mathcal E).
\end{tikzcd}
\]
\end{Proposition}
\begin{proof}
Both lower-right objects represent the functor sending a graded $I$-algebra sheaf $\mathcal B$ on $U$ to graded module maps $j^{-1}\mathcal E\to\mathcal B$. Yoneda gives the unique displayed isomorphism.
\end{proof}

\subsection{Presheaves, stalks and sheafification}
A presheaf of graded algebras satisfying $I$ is a contravariant functor
\[
\mathcal A:\operatorname{Op}(X)^{\op}\longrightarrow G\text{-}\mathbf{Alg}_I.
\]
The use of the opposite category is essential: an inclusion $V\subseteq U$ gives a restriction map $\mathcal A(U)\to\mathcal A(V)$.  The stalk is the filtered colimit
\[
\mathcal A_x=\varinjlim_{x\in U}\mathcal A(U).
\]
By Proposition~\ref{prop:equational-closure}, every stalk satisfies $I$.

\begin{Proposition}\label{prop:sheafification-identities}
Let $\mathcal A$ be a presheaf of graded $F$-algebras satisfying $I$ on every open set.  Its sheafification in graded $F$-modules carries a unique sheaf algebra structure extending that of $\mathcal A$, and the resulting sheaf satisfies $I$.  The canonical map induces isomorphisms on stalks.
\end{Proposition}
\begin{proof}
Addition, multiplication, the unit and the homogeneous projections are morphisms of presheaves and therefore sheafify.  Associativity and the unit identities can be checked on stalks.  The same is true for every polynomial identity: a local evaluation has zero germ at every point, hence is the zero section.  The standard sheafification map is a stalkwise isomorphism.
\end{proof}

Thus objectwise algebraic constructions must be sheafified unless the sheaf axiom has been verified directly.

\subsection{The category of PI-sheaves and the explicit stalk construction}
Let $G\text{-}\mathbf{Alg}_I$ denote the category of graded algebras satisfying $I$.  We write
\[
\operatorname{PSh}_I^G(X)=\Fun(\operatorname{Op}(X)^{\op},G\text{-}\mathbf{Alg}_I)
\]
and let $\operatorname{Sh}_I^G(X)$ be its full subcategory of sheaves.  Limits are computed on open sets.  Colimits are first computed in presheaves and then sheafified.  In particular, the inclusion of sheaves into presheaves has a left adjoint.

The stalk also admits the following explicit description.  An element of $\mathcal A_x$ is an equivalence class $[U,s]_x$, where $x\in U$ and $s\in\mathcal A(U)$; two pairs are equivalent if their restrictions agree on some neighbourhood of $x$ contained in the intersection.  Addition, multiplication and the grading are represented after restriction to a common neighbourhood.  If $f\in I$ and $[U_i,s_i]_x$ are homogeneous germs of the prescribed degrees, choose $V\subseteq\bigcap_iU_i$ containing $x$.  Then
\[
f([s_1]_x,\ldots,[s_m]_x)=
[f(s_1|_V,\ldots,s_m|_V)]_x=0.
\]
Thus the stalk remains in the same equational variety.  This observation gives the local algebraic compatibility between the sheaf and moduli viewpoints.

\begin{Theorem}\label{thm:PI-sheafification}
The sheafification functor
$
a:\operatorname{PSh}_I^G(X)\longrightarrow\operatorname{Sh}_I^G(X)
$
is left adjoint to the inclusion.  It preserves finite products and induces isomorphisms on stalks.  If $\mathcal A$ is already a sheaf after forgetting multiplication, the presheaf multiplication sheafifies uniquely and no further quotient is required.
\end{Theorem}
\begin{proof}
Apply ordinary sheafification to each homogeneous component.  The operations of the algebra are morphisms of presheaves, and the finite-limit preservation of sheafification supplies the induced operations.  All axioms and identities are stalkwise equalities.  The ordinary universal property gives the adjunction, and the stalk statement is standard.
\end{proof}

\subsection{Centrally locally ringed PI-spaces}
\begin{Definition}\label{def:central-PI-space}
A centrally ringed graded PI-space of type $I$ is a triple $(X,\mathcal O_X,\mathcal A)$ such that $(X,\mathcal O_X)$ is a commutative ringed space, $\mathcal A$ is a sheaf of unital graded $\mathcal O_X$-algebras satisfying $I$, and the structural morphism has image in the centre $Z(\mathcal A)$.  It is centrally locally ringed if $(X,\mathcal O_X)$ is locally ringed.  It is finite if $\mathcal A$ is locally finitely generated as an $\mathcal O_X$-module.
\end{Definition}

The definition imposes localness on the centre.  For example, $M_n(\mathcal O_{X,x})$ is not a local ring when $n>1$, but its centre is the local ring $\mathcal O_{X,x}$.  Azumaya and Clifford sheaves therefore define centrally locally ringed PI-spaces.

\begin{Proposition}\label{prop:matrix-stalk-semilocal}
Let $(X,\mathcal O_X)$ be locally ringed and let $A$ be a finite-dimensional $F$-algebra.  Then $(X,\mathcal O_X,A\otimes_F\mathcal O_X)$ is a finite centrally locally ringed PI-space.  If $A=M_n(F)$, then
\[
J(M_n(\mathcal O_{X,x}))=M_n(\mathfrak m_x),\qquad
M_n(\mathcal O_{X,x})/J\cong M_n(k(x)).
\]
\end{Proposition}
\begin{proof}
The first statement is Proposition~\ref{prop:tensor-constant-finite}.  The radical formula for matrix rings is standard and the quotient formula follows by reduction modulo the maximal ideal.
\end{proof}

\subsection{Recovery of super- and colour-geometric structure sheaves}
Let $\varepsilon:G\times G\to F^*$ be a bicharacter and let $I_\varepsilon$ be the graded $T$-ideal generated by $x^gy^h-\varepsilon(g,h)y^hx^g$.
\begin{Proposition}\label{prop:colour-sheaf-equivalence}
For every ringed space $(X,\mathcal O_X)$, the full category $\operatorname{Sh}_{I_\varepsilon}^G(X)$ is canonically equal to the category of sheaves of unital $\varepsilon$-commutative $G$-graded $\mathcal O_X$-algebras.  In particular, the cases
$
(G,\varepsilon)=(\mathbb Z_2,(-1)^{ij})
 $and $\quad
G=\mathbb Z_2^n
$
recover, respectively, the algebraic structure sheaves of supergeometry and colour supergeometry.
\end{Proposition}
\begin{proof}
A degree-preserving evaluation of the generators of $I_\varepsilon$ is exactly the relation $ab=\varepsilon(g,h)ba$ for homogeneous sections $a\in\mathcal A_g$ and $b\in\mathcal A_h$.  Since $I_\varepsilon$ is the graded $T$-ideal generated by these elements, vanishing of the generators is equivalent to vanishing of every element of $I_\varepsilon$.
\end{proof}

Proposition~\ref{prop:colour-sheaf-equivalence} does not encode the local-model conditions of a supermanifold, which concern the degree-zero body, the nilpotent or adic ideal, smoothness and completeness.  They do not require the total graded algebra to be a local ring in the ordinary noncommutative sense.

\begin{Corollary}\label{cor:split-super-classifying-map}
Let $\mathcal V$ be a locally free odd vector bundle of rank $q$ on a scheme $X$.  The split superalgebra
\[
\mathcal A=\bigwedge_{\mathcal O_X}\mathcal V^\vee
\]
is finite locally free of superrank $(2^{q-1},2^{q-1})$ for $q>0$, satisfies $I_\varepsilon$, and determines an object of
\[
\PI^{\mathbb Z_2}_{I_\varepsilon,(2^{q-1},2^{q-1})}(X).
\]
For $q=0$ it determines the purely even rank-one point.
\end{Corollary}
\begin{proof}
Locally trivialise $\mathcal V$.  The exterior algebra then has the stated even and odd ranks and satisfies the supercommutativity identities.  Theorem~\ref{thm:PI-stack-classification} supplies the corresponding classifying object.
\end{proof}

\subsection{The sheaf and stack viewpoints}
\begin{Proposition}\label{prop:sheaf-stack-interface}
Let $S$ be an $F$-scheme.  The groupoid of sheaves $\mathcal B=\bigoplus_g\mathcal B_g$ of unital graded $\mathcal O_S$-algebras satisfying $I$, with each $\mathcal B_g$ locally free of rank $d_g$, is naturally equivalent to the groupoid of morphisms
\[
S\longrightarrow\PI^G_{I,\mathbf d}.
\]
Under this equivalence, pullback of families corresponds to composition of classifying morphisms.
\end{Proposition}
\begin{proof}
This is the functorial form of Theorem~\ref{thm:PI-stack-classification}.
\end{proof}
\section{The moduli stack of graded PI-algebra laws}\label{sec:pi-moduli}

Throughout this section, $I$ is a finitely generated graded $T$-ideal; by Theorem~\ref{thm:graded-specht}, this includes $I=\idG(A)$ for every graded PI-algebra $A$.  Let $\mathbf d=(d_g)_{g\in G}$ be a dimension vector with $d_e\geq1$.  Choose graded vector spaces $V_g$ with $\dim_FV_g=d_g$ and a distinguished vector $\mathbf 1\in V_e$.  Put $V=\bigoplus_gV_g$ and
\[
H_{\mathbf d}:=\{(h_g)_g\in\prod_{g\in G}\GL(V_g):h_e(\mathbf1)=\mathbf1\}.
\]

A graded multiplication law on $(V,\mathbf1)$ is a family $\mu_{g,h}:V_g\otimes V_h\longrightarrow V_{gh}$.
The coefficients of $\mu$ form an affine space
\[
\mathbb M^G_{\mathbf d}:=
\bigoplus_{g,h\in G}\operatorname{Hom}_F(V_g\otimes V_h,V_{gh}).
\]
Associativity and the requirement that $\mathbf1$ be a two-sided unit are polynomial equations in these coefficients.

Choose multilinear graded $T$-generators $f_1,\ldots,f_r$ of $I$.  For every $f_j$ and every degree-compatible tuple of basis elements of $V$, expand $f_j$ using $\mu$.  The vanishing of all coordinates of all such evaluations is again a finite system of polynomial equations.

\begin{Definition}\label{def:coordinate-ring-law-scheme}
Let
\[
P_{\mathbf d}=\operatorname{Sym}_F\bigl((\mathbb M^G_{\mathbf d})^\vee\bigr)
=F[c^{g,h,k}_{i,j}]
\]
be the polynomial coordinate ring of the affine space of graded multiplications.  Let $\mathfrak a_{\mathrm{ass}}$, $\mathfrak a_{\mathrm{unit}}$ and $\mathfrak a_I$ be, respectively, the ideals generated by the coefficients of the associators, the two unit equations and the evaluations of a finite multihomogeneous generating system of $I$.  Put
\[
R^G_{I,\mathbf d}
:=P_{\mathbf d}/(\mathfrak a_{\mathrm{ass}}+\mathfrak a_{\mathrm{unit}}+\mathfrak a_I).
\]
\end{Definition}

\begin{Proposition}\label{prop:universal-law-coordinate-ring}
The free graded $R^G_{I,\mathbf d}$-module
\[
\mathcal U^{\mathrm{fr}}_{I,\mathbf d}
=V\otimes_F R^G_{I,\mathbf d}
\]
carries a universal unital graded multiplication $\mu^{\mathrm{univ}}$ satisfying $I$.  For every commutative $F$-algebra $R$, evaluation of the coordinate functions induces a natural bijection
\[
\operatorname{Hom}_{F\text{-}\mathrm{Alg}}(R^G_{I,\mathbf d},R)
\cong
\left\{\begin{array}{c}
\text{unital graded $R$-algebra laws on $V\otimes_F R$}\\
\text{with fixed unit $\mathbf1$ and satisfying $I$}
\end{array}\right\}.
\]
\end{Proposition}
\begin{proof}
The structure constants $c^{g,h,k}_{i,j}$ define the universal bilinear law on $V\otimes_FP_{\mathbf d}$.  Passing to $R^G_{I,\mathbf d}$ imposes precisely the associativity, unit and identity equations.  A homomorphism $R^G_{I,\mathbf d}\to R$ is therefore exactly a compatible specialisation of all structure constants, and the construction is natural in $R$.
\end{proof}

\begin{Definition}
The affine scheme of unital graded algebra laws satisfying $I$ is the closed subscheme
\[
\AlgLaws^G_{\mathbf d}(I)\subseteq\mathbb M^G_{\mathbf d}
\]
defined by the associativity, unit and PI equations.  The group $H_{\mathbf d}$ acts by change of basis.  Using the standard quotient-stack convention \cite{VistoliNotes}, we define
\[
\PI^G_{I,\mathbf d}:=
[\AlgLaws^G_{\mathbf d}(I)/H_{\mathbf d}].
\]
If $I=\idG(A)$, we write $\PI^G_{A,\mathbf d}$.
\end{Definition}

\begin{Proposition}\label{prop:finite-type-law-scheme}
The scheme $\AlgLaws^G_{\mathbf d}(I)$ is affine of finite type over $F$ and is independent of the chosen finite $T_G$-generating system of $I$.
\end{Proposition}

\begin{proof}
The ambient space is finite-dimensional affine space.  Associativity and unit give finitely many equations.  By Theorem~\ref{thm:graded-specht}, $I$ has finitely many $T_G$-generators, and each generator yields finitely many coefficient equations because $V$ is finite-dimensional.  A law satisfies one finite generating system if and only if it satisfies the whole $T$-ideal, hence if and only if it satisfies any other generating system.  The defining functor and therefore the closed subscheme are independent of the choice.
\end{proof}

\begin{Theorem}\label{thm:PI-stack-algebraicity}
The quotient $\PI^G_{I,\mathbf d}$ is an algebraic stack of finite type over $F$.  More precisely:
\begin{enumerate}
\item the canonical morphism
\[
\pi:\AlgLaws^G_{\mathbf d}(I)\longrightarrow\PI^G_{I,\mathbf d}
\]
is a smooth surjective atlas of relative dimension
\[
\dim H_{\mathbf d}=\sum_{g\in G}d_g^2-d_e;
\]
\item the diagonal
\[
\Delta:\PI^G_{I,\mathbf d}\longrightarrow
\PI^G_{I,\mathbf d}\times_F\PI^G_{I,\mathbf d}
\]
is representable by affine schemes of finite presentation;
\item for every geometric point $[B]:\Spec k\to\PI^G_{I,\mathbf d}$, the stabiliser group scheme is $\operatorname{Stab}_{[B]}\simeq\Aut_{G,1}(B)$, the affine group scheme of degree-preserving unital algebra automorphisms of $B$, and $\Lie\Aut_{G,1}(B)\simeq\Der_G(B)$.
\end{enumerate}
\end{Theorem}
\begin{proof}
Write $X_I=\AlgLaws^G_{\mathbf d}(I)$ and $H=H_{\mathbf d}$.  Choose a
splitting $V_e=F\mathbf1\oplus W$.  Relative to this splitting, an element
of $\GL(V_e)$ fixing $\mathbf1$ has a unique block form
\[
\begin{pmatrix}1&\ell\\0&A\end{pmatrix},
\qquad \ell\in\operatorname{Hom}_F(W,F),\quad A\in\GL(W).
\]
Thus
\[
H\cong
\bigl(\operatorname{Hom}_F(W,F)\rtimes\GL(W)\bigr)
\times\prod_{g\ne e}\GL(V_g).
\]
In particular, $H$ is a smooth affine group scheme of finite type and
\[
\dim H=(d_e-1)+(d_e-1)^2+\sum_{g\ne e}d_g^2
=\sum_{g\in G}d_g^2-d_e.
\]

The action groupoid presenting the quotient stack is
\[
H\times_FX_I\;\substack{\xrightarrow{\ a\ }\\[-.4ex]
\xrightarrow[\ \mathrm{pr}_2\ ]{}}\;X_I.
\]
Both arrows are smooth and of finite presentation because $H$ is smooth.
The general quotient-stack theorem therefore applies
\cite{VistoliNotes,StacksQuotient}: $[X_I/H]$ is an algebraic stack, and
$\pi:X_I\to[X_I/H]$ is a smooth surjective morphism.  After any base change
$T\to[X_I/H]$, the morphism $X_I\times_{[X_I/H]}T\to T$ is the
$H_T$-torsor attached to the object over $T$; hence it is smooth of relative
dimension $\dim H$.  Since $X_I$ is of finite type over $F$, so is the
quotient stack.  This proves (1).

We next prove (2) directly.  Let $S$ be an $F$-scheme and let
$\mathcal B,\mathcal C$ be two $S$-objects of the quotient stack, viewed by
Theorem~\ref{thm:PI-stack-classification} as graded algebra bundles.  The
functor
\[
T\longmapsto
\operatorname{Isom}_{G,1;\,\mathcal O_T\text{-alg}}
(\mathcal B_T,\mathcal C_T)
\]
is a closed subfunctor of the functor of degree-preserving vector-bundle
isomorphisms carrying $1_{\mathcal B}$ to $1_{\mathcal C}$.  The latter is
an affine $S$-scheme of finite presentation: locally on $S$ it is the
translate of the affine group $H_S$.  In homogeneous local frames, the
additional conditions are the finitely many equations
\[
\phi\bigl(\mu_{\mathcal B}(e_i^g,e_j^h)\bigr)
=\mu_{\mathcal C}\bigl(\phi(e_i^g),\phi(e_j^h)\bigr)
\]
for all basis vectors.  They are polynomial equations in the matrix
coefficients of $\phi$.  Hence the isomorphism functor is represented by an
affine $S$-scheme of finite presentation.  Since these isomorphism schemes
are precisely the base changes of the diagonal, the diagonal is
representable, affine, and of finite presentation.

Finally let $B$ be a geometric point.  The fibre of the diagonal over
$([B],[B])$ represents degree-preserving unital algebra automorphisms of
$B$, so it is $\Aut_{G,1}(B)$.  To compute its Lie algebra, let
$D=k[\varepsilon]/(\varepsilon^2)$.  An element in the tangent space at the
identity has the form $g=\id+\varepsilon\delta$, where $\delta$ is a
degree-preserving linear endomorphism with $\delta(1)=0$.  The equality
$g(xy)=g(x)g(y)$ modulo $\varepsilon^2$ is equivalent to
\[
\delta(xy)=\delta(x)y+x\delta(y).
\]
Thus the tangent space at the identity, with its commutator bracket, is
$\Der_G(B)$.  This proves (3).
\end{proof}

\begin{Corollary}\label{cor:PI-DM-locus}
There is an open substack
$
\PI^{G,\mathrm{DM}}_{I,\mathbf d}\subseteq\PI^G_{I,\mathbf d}
$
whose geometric points are the algebras $B$ satisfying
$
\Der_G(B)=0.
$
It is the maximal open substack on which the diagonal is unramified; equivalently, its geometric stabilisers are finite \'{e}tale group schemes.
\end{Corollary}
\begin{proof}
Every algebraic stack has a largest open Deligne--Mumford substack, and a
geometric point belongs to it exactly when its stabiliser is unramified
\cite{StacksDMLocus}.  By Theorem~\ref{thm:PI-stack-algebraicity}, the
stabiliser at $[B]$ is the affine group scheme of finite type
$\Aut_{G,1}(B)$ and its Lie algebra is $\Der_G(B)$.  Over a field of
characteristic zero, a finite-type group scheme is smooth.  For a smooth
group scheme, unramifiedness is equivalent to relative dimension zero,
which is equivalent to the vanishing of its tangent space at the identity.
Therefore
\[
[B]\text{ lies in the Deligne--Mumford locus}
\quad\Longleftrightarrow\quad
\Der_G(B)=0.
\]
On this locus the stabilisers are affine, smooth, of finite type and zero-dimensional; hence they are finite \'{e}tale.  Maximality follows from the defining property of the Deligne--Mumford locus.
\end{proof}

\begin{Corollary}\label{cor:PI-stack-local-dimension}
Let $k/F$ be algebraically closed and let $B$ be a $k$-point of
$\AlgLaws^G_{\mathbf d}(I)$.  The stack $\PI^G_{I,\mathbf d}$ is smooth at $[B]$ if and only if the law scheme is smooth at $B$.  In that case
\[
\dim_{[B]}\PI^G_{I,\mathbf d}
=
\dim_B\AlgLaws^G_{\mathbf d}(I)-\dim H_{\mathbf d}
=
\dim Z^2_{G,I,\mathrm n}(B,B)-\dim C^1_{G,\mathrm n}(B,B),
\]
and equivalently
\[
\dim_{[B]}\PI^G_{I,\mathbf d}
=
\dim H^0\bigl(\mathbb T^{\mathrm{cl}}_{I,[B]}\bigr)
-
\dim H^{-1}\bigl(\mathbb T^{\mathrm{cl}}_{I,[B]}\bigr).
\]
\end{Corollary}
\begin{proof}
Let $X_I=\AlgLaws^G_{\mathbf d}(I)$ and $H=H_{\mathbf d}$.  The atlas
$\pi:X_I\to[X_I/H]$ is smooth and open.  If the stack is smooth over $F$ at
$[B]$, then $X_I$ is smooth over $F$ at $B$ by composition with $\pi$.
Conversely, assume $X_I$ is smooth at $B$.  The smooth locus
$X_I^{\mathrm{sm}}\subseteq X_I$ is open and $H$-invariant, because the
$H$-action is by $F$-scheme automorphisms.  Hence
$[X_I^{\mathrm{sm}}/H]$ is an open substack containing $[B]$, and
$X_I^{\mathrm{sm}}\to[X_I^{\mathrm{sm}}/H]$ is smooth and surjective.
Since $X_I^{\mathrm{sm}}$ is smooth over $F$, smoothness descends through
this atlas.  Thus the stack is smooth at $[B]$, proving the first assertion.

At a smooth point, local dimension equals tangent-space dimension.  The
relative dimension of the atlas is $\dim H$, hence the dimension formula for
a smooth presentation gives
\[
\dim_{[B]}[X_I/H]=\dim_BX_I-\dim H.
\]
Linearising associativity, the unit, the grading and the chosen PI-equations
shows that
\[
T_BX_I=Z^2_{G,I,\mathrm n}(B,B).
\]
The Lie algebra of $H$ consists of degree-preserving endomorphisms of $B$
which vanish on the unit; this is $C^1_{G,\mathrm n}(B,B)$.  Therefore
\[
\dim_{[B]}\PI^G_{I,\mathbf d}
=\dim Z^2_{G,I,\mathrm n}(B,B)
-\dim C^1_{G,\mathrm n}(B,B).
\]
For a finite-dimensional two-term complex $C^{-1}\to C^0$, rank-nullity
gives
\[
\dim H^0-\dim H^{-1}=\dim C^0-\dim C^{-1}.
\]
Applying this to $\mathbb T^{\mathrm{cl}}_{I,[B]}$ proves the final formula.
\end{proof}

\subsection{The fibred moduli category and its descent diagram}
Define a category fibred in groupoids
\[
\mathcal{PI}_{I,\mathbf d}\longrightarrow(\mathbf{Sch}/F)
\]
by letting $\mathcal{PI}_{I,\mathbf d}(S)$ be the groupoid of locally free graded $\mathcal O_S$-algebras of graded rank $\mathbf d$ satisfying $I$. Pullback along $f:T\to S$ is
\[
f^*(\mathcal B,\mu,1)=(f^*\mathcal B,f^*\mu,f^*1).
\]
The classification theorem identifies this fibred category with the quotient stack $[\AlgLaws^G_{\mathbf d}(I)/H_{\mathbf d}]$.
For a family $\mathcal B$ on $S$, let
$
P=\underline{\operatorname{Isom}}_{G,1}
(V\otimes_F\mathcal O_S,\mathcal B)
$
be the sheaf of graded frames carrying $\mathbf1$ to the unit. The descent datum is encoded by
\[
\begin{tikzcd}[column sep=large,row sep=large]
P\times H_{\mathbf d}
\arrow[r,shift left=.8ex,"\mathrm{pr}_1"]
\arrow[r,shift right=.8ex,"a"']
& P \arrow[r,"\pi"] \arrow[d,"\widetilde\mu"]
& S\\
& \AlgLaws^G_{\mathbf d}(I),
\end{tikzcd}
\]
where $a$ is the right action and $\widetilde\mu$ sends a frame to the transported multiplication. Equivariance is
$
\widetilde\mu(p\cdot h)=h^{-1}\cdot\widetilde\mu(p).
$
Conversely, an equivariant map from an $H_{\mathbf d}$-torsor descends the trivial bundle $P\times V$ and its universal multiplication.

\begin{Proposition}\label{prop:classifying-base-change-square}
Let $c_{\mathcal B}:S\to\PI^G_{I,\mathbf d}$ classify $\mathcal B$, and let $f:T\to S$. Then $f^*\mathcal B$ is classified by a $2$-commutative square
\[
\begin{tikzcd}[column sep=large,row sep=large]
T \arrow[r,"c_{f^*\mathcal B}"] \arrow[d,"f"']
& \PI^G_{I,\mathbf d} \arrow[d,equal]\\
S \arrow[r,"c_{\mathcal B}"']
& \PI^G_{I,\mathbf d}.
\end{tikzcd}
\]
The induced isomorphism between the pullbacks of the tautological family is the canonical base-change isomorphism.
\end{Proposition}
\begin{proof}
The frame torsor of $f^*\mathcal B$ is $P\times_ST$, and its equivariant law map is the pullback of $\widetilde\mu$. This is the composite of the classifying morphism with $f$.
\end{proof}

\begin{Example}\label{ex:two-dimensional-chart}
Take the trivial grading and $V=F\mathbf1\oplus Fx$.  Every unital algebra law on this module is uniquely $x^2=a\mathbf1+bx$, with $(a,b)\in\mathbb A^2_F$, and is the associative commutative quotient $B_{a,b}=F[t]/(t^2-bt-a)$.  Thus the rank-two commutative law scheme is $\mathbb A^2_F$.  A general unit-preserving change of basis is $x'=r\mathbf1+\lambda x$, where $r\in F$ and $\lambda\in F^*$, and it sends $(a,b)$ to $(a',b')=(\lambda^2a-r^2-\lambda br,\,2r+\lambda b)$.  Consequently the discriminant $\Delta=b^2+4a$ transforms as $\Delta'=\lambda^2\Delta$.  Hence $D(\Delta)$ descends to the étale-semisimple locus.  Over an algebraically closed field of characteristic different from $2$, its fibres are $F\times F$, while $V(\Delta)$ is the dual-number degeneration.  This chart exhibits the equations, orbits, stabilisers and strata explicitly.
\end{Example}

\begin{Theorem}\label{thm:PI-stack-classification}
For every $F$-scheme $S$, the groupoid $\PI^G_{I,\mathbf d}(S)$ is naturally equivalent to the groupoid of pairs $(\mathcal B,\iota)$, where
\begin{enumerate}
\item $\mathcal B=\bigoplus_{g\in G}\mathcal B_g$ is a unital graded $\mathcal O_S$-algebra;
\item every $\mathcal B_g$ is locally free of rank $d_g$;
\item $\iota:\mathcal O_S\to\mathcal B_e$ is the unit section;
\item every polynomial in $I$ vanishes under every degree-preserving local evaluation in $\mathcal B$.
\end{enumerate}
Morphisms are graded $\mathcal O_S$-algebra isomorphisms.
\end{Theorem}

\begin{proof}
Set $X_I=\AlgLaws^G_{\mathbf d}(I)$ and $H=H_{\mathbf d}$.  We construct
quasi-inverse functors between the two groupoids.
An object of $[X_I/H](S)$ is an $H$-torsor $p:P\to S$ for the fppf topology
together with an equivariant morphism $u:P\to X_I$.  For each $g\in G$ set
\[
\mathcal B_g=P\times^H V_g.
\]
These are locally free $\mathcal O_S$-modules of ranks $d_g$.  The universal
multiplication on $V\otimes\mathcal O_{X_I}$ pulls back to an
$H$-equivariant multiplication on $P\times V$; effective descent for
quasi-coherent modules and their morphisms produces maps
\[
\mu_{g,h}:\mathcal B_g\otimes\mathcal B_h\longrightarrow\mathcal B_{gh}.
\]
The fixed vector $\mathbf1\in V_e$ descends to a section
$1_{\mathcal B}\in\mathcal B_e$.  Associativity, the two unit identities,
and every evaluation of the generators of $I$ hold after pullback to the
faithfully flat cover $P\to S$, because $u$ factors through $X_I$.  Equality
of morphisms of quasi-coherent sheaves descends faithfully, so the same
identities hold on $S$.  Hence $\mathcal B=\bigoplus_g\mathcal B_g$ is a
family of the required kind.

Conversely, let $\mathcal B$ be such a family.  Define the fppf sheaf
$
P=\underline{\operatorname{Isom}}_{G,1}
(V\otimes_F\mathcal O_S,\mathcal B).$ 
First, the frame torsor is locally nonempty.  Each $\mathcal B_g$ is locally free.  For the
neutral component, the unit is nonzero in every geometric fibre: otherwise
the fibre algebra would be zero, contradicting $d_e\ge1$.  In a local frame
of $\mathcal B_e$, the coordinates of $1_{\mathcal B}$ therefore generate
the unit ideal.  After shrinking, one coordinate is invertible, so elementary
column operations extend $1_{\mathcal B}$ to a basis.  Thus, Zariski-locally
on $S$, there are homogeneous frames carrying $\mathbf1$ to
$1_{\mathcal B}$.  It follows that $P$ is an $H$-torsor (indeed already for
the Zariski topology after the individual bundles have been trivialised).

A frame $p\in P(T)$ transports the multiplication of $\mathcal B_T$ to a
law on $V\otimes_F\mathcal O_T$.  This gives a morphism
$u:P\to\mathbb M^G_{\mathbf d}$.  The transported law is associative,
unital, graded, and satisfies $I$, so $u$ factors through $X_I$.  If the
right action on frames is $p\cdot h=p\circ h$, then transport gives
$u(p\cdot h)=h^{-1}\cdot u(p)$, which is precisely the equivariance
convention for the quotient stack.
The two constructions are quasi-inverse.  Starting with $(P,u)$, the map
\[
P\times V\longrightarrow p^*(P\times^H V),\qquad
(p,v)\longmapsto(p,[p,v]),
\]
is the tautological equivariant trivialisation; it identifies the recovered
frame torsor and law map with $(P,u)$.  Starting with $\mathcal B$, the map
\[
P\times^H V\longrightarrow\mathcal B,\qquad [p,v]\longmapsto p(v),
\]
is well defined, is an isomorphism after pulling back to $P$, and respects
the multiplication and unit; faithful flat descent makes it an algebra
isomorphism on $S$.

Finally, an isomorphism of algebra families carries frames to frames and
therefore induces a morphism of the corresponding torsors compatible with
the law maps.  Conversely, a morphism of torsor data descends to an algebra
isomorphism of associated bundles.  These operations are inverse on
morphisms as well as on objects and commute with arbitrary base change in
$S$.  Hence they define the claimed natural equivalence of groupoids.
\end{proof}

\begin{Definition}
A finite-rank PI-geometry of type $(I,\mathbf d)$ over $S$ is a morphism
\[
S\longrightarrow\PI^G_{I,\mathbf d}.
\]
By Theorem~\ref{thm:PI-stack-classification}, this is equivalently a locally free family of graded algebras satisfying $I$.
\end{Definition}

\begin{Theorem}\label{thm:closed-PI-immersion}
If $I\subseteq J$ are finitely generated graded $T$-ideals, then the additional identities in $J$ define an $H_{\mathbf d}$-invariant closed immersion
\[
\AlgLaws^G_{\mathbf d}(J)\hookrightarrow\AlgLaws^G_{\mathbf d}(I),
\]
and hence a canonical closed immersion
$
\PI^G_{J,\mathbf d}\hookrightarrow\PI^G_{I,\mathbf d}.
$
These immersions are functorial for chains of $T_G$-ideals.
\end{Theorem}

\begin{proof}
Use the common polynomial coordinate ring $P_{\mathbf d}$ of all graded
multiplications.  Let $\mathfrak a_I$ and $\mathfrak a_J$ denote the ideals
generated by associativity, the unit equations, and the coefficient
functions obtained from evaluations of $I$ and $J$, respectively.  Since
$I\subseteq J$, every coefficient equation imposed by $I$ is also imposed
by $J$, and therefore
\[
\mathfrak a_I\subseteq\mathfrak a_J.
\]
The quotient homomorphism
$P_{\mathbf d}/\mathfrak a_I\twoheadrightarrow
P_{\mathbf d}/\mathfrak a_J$ represents a closed immersion
$
\AlgLaws^G_{\mathbf d}(J)\hookrightarrow
\AlgLaws^G_{\mathbf d}(I).
$
It is $H_{\mathbf d}$-equivariant because polynomial identities are
preserved by graded change of basis.
The induced morphism of quotient stacks is representable by this closed
immersion.  Indeed, after base change by the smooth atlas
$\AlgLaws^G_{\mathbf d}(I)\to\PI^G_{I,\mathbf d}$, one obtains exactly
\[
\AlgLaws^G_{\mathbf d}(J)\to\AlgLaws^G_{\mathbf d}(I).
\]
Closed immersions are local on the target for the smooth topology, so
$\PI^G_{J,\mathbf d}\to\PI^G_{I,\mathbf d}$ is a closed immersion.  For a
chain $I\subseteq J\subseteq K$, all three maps are induced by the nested
coordinate ideals
$\mathfrak a_I\subseteq\mathfrak a_J\subseteq\mathfrak a_K$, which proves
functoriality.
\end{proof}

\begin{Proposition}\label{prop:comparison-2-cartesian}
For $I\subseteq J$, the diagram
\[
\begin{tikzcd}[column sep=large,row sep=large]
\AlgLaws^G_{\mathbf d}(J) \arrow[r,hook] \arrow[d]
& \AlgLaws^G_{\mathbf d}(I) \arrow[d]\\
\PI^G_{J,\mathbf d} \arrow[r,hook]
& \PI^G_{I,\mathbf d}
\end{tikzcd}
\]
is $2$-Cartesian. The pullback of the tautological $I$-algebra to $\PI^G_{J,\mathbf d}$ is canonically the tautological $J$-algebra.
\end{Proposition}
\begin{proof}
The vertical maps are the universal $H_{\mathbf d}$-torsors. Base change of the atlas $\AlgLaws^G_{\mathbf d}(I)\to[\AlgLaws^G_{\mathbf d}(I)/H_{\mathbf d}]$ along $[\AlgLaws^G_{\mathbf d}(J)/H_{\mathbf d}]$ is the torsor $\AlgLaws^G_{\mathbf d}(J)\to[\AlgLaws^G_{\mathbf d}(J)/H_{\mathbf d}]$. Both tautological families descend from $V\otimes\mathcal O_{\AlgLaws^G_{\mathbf d}(J)}$ with the same law.
\end{proof}

\subsection{Equations in structure constants}
Choose homogeneous bases $e^g_1,\ldots,e^g_{d_g}$ and write
\[
\mu(e^g_i,e^h_j)=\sum_{k=1}^{d_{gh}}c^{g,h,k}_{i,j}e^{gh}_k.
\]
Associativity is the finite collection of quadratic equations
\[
\sum_{\ell}c^{g,h,\ell}_{i,j}c^{gh,r,m}_{\ell,k}
=
\sum_{\ell}c^{h,r,\ell}_{j,k}c^{g,hr,m}_{i,\ell}.
\]
If the distinguished unit is $e^e_1$, the unit equations are
$
c^{e,g,k}_{1,j}=\delta_{jk}, $ and $ c^{g,e,k}_{i,1}=\delta_{ik}.
$
A multilinear identity $f(x_1^{g_1},\ldots,x_m^{g_m})$ becomes, after evaluation on basis vectors, finitely many polynomial equations in the $c^{g,h,k}_{i,j}$.  This presentation proves that the scheme is of finite type and retains scheme-theoretic multiplicities that the set of isomorphism classes would discard.

\begin{Example}
For the trivial grading and rank $d$, the ambient associative-law scheme is cut out by associativity and the unit.  Adding $[x,y]$ imposes the linear equations $c^k_{ij}=c^k_{ji}$.  The resulting quotient stack is the moduli stack of commutative algebra structures on rank-$d$ vector bundles.  Its nonreduced structure records infinitesimal commutative deformations.
\end{Example}

\begin{Example}
For $G=\mathbb Z_2$ and superdimension $(p,q)$, the multiplication coordinates occur only in the degree-compatible blocks.  The equations
$
\mu(x,y)=(-1)^{|x||y|}\mu(y,x)
$
are linear in the structure constants, while associativity remains quadratic.  The quotient by $\GL_p\times\GL_q$ with the unit fixed is the stack of finite locally free supercommutative algebra laws of superrank $(p,q)$.
\end{Example}

\subsection{The tautological family and descent}
Unlike a coarse orbit space, the quotient stack retains stabilisers and supports a tautological family.  The atlas $\AlgLaws^G_{\mathbf d}(I)\to\PI^G_{I,\mathbf d}$ carries the trivial graded vector bundle $V\otimes_F\mathcal O$ with its universal multiplication; $H_{\mathbf d}$-equivariance descends it to a graded algebra $\mathcal U_{I,\mathbf d}$ on the stack.

\begin{Proposition}\label{prop:tautological-family}
For every morphism $f:S\to\PI^G_{I,\mathbf d}$, the associated family is canonically isomorphic to $f^*\mathcal U_{I,\mathbf d}$.  Conversely, assigning a classifying morphism to a family defines an equivalence of groupoids as in Theorem~\ref{thm:PI-stack-classification}; in particular, the classifying morphism is determined up to $2$-isomorphism, and its $2$-automorphisms are exactly the automorphisms of the family.
\end{Proposition}
\begin{proof}
This is the descent construction in Theorem~\ref{thm:PI-stack-classification}: after pulling back to the frame torsor, the family becomes the trivial bundle equipped with the corresponding law.  Descent identifies it with the pullback of the tautological object.
\end{proof}

\subsection{Relation with representation schemes and generic matrices}
The law stack parametrises algebra structures on a fixed-rank bundle, whereas $\operatorname{Rep}_n(C)$ parametrises homomorphisms from a fixed algebra $C$ to $M_n(-)$.  Universal matrix representations of relatively free algebras produce generic matrix algebras; conversely, the regular representation sends a family of laws to endomorphism algebras.  Thus representation schemes record modules over a fixed algebra, while the quotient stack records deformations of the coefficient algebra.  On the Azumaya locus, Morita theory reduces the regular representation to the centre.

\begin{Proposition}\label{prop:regular-representation-family}
Let $\mathcal B$ be a finite locally free algebra of rank $d$ on $S$.  Left multiplication defines a faithful morphism of $\mathcal O_S$-algebras
$
\lambda:\mathcal B\longrightarrow\cEnd_{\mathcal O_S}(\mathcal B).
$
It is compatible with arbitrary base change and with the classifying morphism $S\to\PI^G_{I,\mathbf d}$.
\end{Proposition}
\begin{proof}
The map sends $b$ to $(x\mapsto bx)$.  Faithfulness follows from evaluating at the unit.  Both multiplication and the endomorphism bundle commute with base change for finite locally free modules.
\end{proof}

\section{PI-restricted deformations and infinitesimal K-theory}\label{sec:pi-deformations}

Before describing the tangent complex of the moduli stack, we recall the universal differential calculus in the form compatible with sheaves and bimodules.  Let $B$ be an $F$-algebra and let $B^e=B\otimes_F B^{\op}$.

\begin{Definition}
The bimodule of universal noncommutative one-forms is
\[
\Omega^1_{B/F}:=\ker\bigl(\mu:B\otimes_F B\to B\bigr),
\qquad
db=1\otimes b-b\otimes1.
\]
For a sheaf $\mathcal B$ of $F$-algebras, the same definition is made in sheaves and gives a $\mathcal B$-bimodule $\Omega^1_{\mathcal B/F}$.
\end{Definition}

\begin{Proposition}\label{prop:universal-differentials}
For every $B$-bimodule $M$, composition with $d$ gives a natural isomorphism
\[
\operatorname{Hom}_{B^e}(\Omega^1_{B/F},M)\xrightarrow{\sim}\Der_F(B,M).
\]
The analogous statement holds for sheaves, with internal sheaf Hom.
\end{Proposition}
\begin{proof}
Let $u:\Omega^1_{B/F}\to M$ be a $B$-bimodule map and define
$\delta_u(b)=u(db)$.  Since
$
d(ab)=a\,db+(da)\,b,
$
we have
$\delta_u(ab)=a\delta_u(b)+\delta_u(a)b$; thus $\delta_u$ is an
$F$-derivation.
Conversely, let $\delta:B\to M$ be an $F$-derivation.  For
$\xi=\sum_i a_i\otimes b_i\in\Omega^1_{B/F}$, so that
$\sum_i a_ib_i=0$, put
$
u_\delta(\xi)=\sum_i a_i\delta(b_i).
$
This is left $B$-linear.  For $c\in B$,
\[
u_\delta(\xi c)=\sum_i a_i\delta(b_ic)=\sum_i a_i\delta(b_i)c+\left(\sum_i a_ib_i\right)\delta(c)=u_\delta(\xi)c.
\]
so it is also right $B$-linear.  Moreover
$u_\delta(db)=\delta(b)$.  Conversely, every
$\xi=\sum_i a_i\otimes b_i\in\ker\mu$ satisfies
\[
\xi=\sum_i a_i(1\otimes b_i-b_i\otimes1)=\sum_i a_i\,db_i,
\]
so a bimodule map is determined by its values on the elements $db$.
Therefore $u\mapsto u\circ d$ and $\delta\mapsto u_\delta$ are inverse.
The construction is natural in $M$.
For sheaves the same formulas are applied on each open set.  They commute
with restriction, and the verification is local, giving the internal sheaf
Hom statement.
\end{proof}

Thus the tangent Lie algebra of an algebra is $\Der_F(B,B)$, or equivalently $\operatorname{Hom}_{B^e}(\Omega^1_{B/F},B)$; it is not a space of derivations \emph{from} $\Omega^1_{B/F}$ treated as an algebra.  Inner derivations yield the exact sequence
\[
0\longrightarrow Z(B)\longrightarrow B\xrightarrow{\operatorname{ad}}\Der_F(B,B)
\longrightarrow HH^1_F(B,B)\longrightarrow0.
\]
On the Azumaya locus this sequence governs derivations, whereas first-order deformations of the multiplication are controlled in degree two.

\subsection{Relative differentials of the PI-envelope}
Put $\mathcal T=\Tens_X^G(\mathcal E)$, $\mathcal J=\mathcal J_I(\mathcal E)$ and $\mathcal B=\mathcal T/\mathcal J=\FreePI^G_{I,X}(\mathcal E)$.
\begin{Proposition}\label{prop:PI-conormal-sequence}
There is a natural right exact sequence of $\mathcal B$-bimodules
\[
\mathcal J/\mathcal J^2
\xrightarrow{\partial}
\mathcal B\otimes_{\mathcal T}
\Omega^1_{\mathcal T/\mathcal O_X}
\otimes_{\mathcal T}\mathcal B
\longrightarrow
\Omega^1_{\mathcal B/\mathcal O_X}
\longrightarrow0,
\]
where $\partial([r])=1\otimes dr\otimes1$.  Consequently, for every $\mathcal B$-bimodule $\mathcal M$, restriction induces an exact sequence
\[
0\longrightarrow\Der_{\mathcal O_X}(\mathcal B,\mathcal M)
\longrightarrow\Der_{\mathcal O_X}(\mathcal T,\mathcal M)
\longrightarrow
\mathcal Hom_{\mathcal B^e}(\mathcal J/\mathcal J^2,\mathcal M).
\]
The final arrow sends a derivation to its values on the defining PI-relations.
\end{Proposition}
\begin{proof}
Write $\mathcal B=\mathcal T/\mathcal J$.  The morphism
\[
\partial:\mathcal J/\mathcal J^2\longrightarrow
\mathcal B\otimes_{\mathcal T}\Omega^1_{\mathcal T/\mathcal O_X}
\otimes_{\mathcal T}\mathcal B,
\qquad [r]\longmapsto1\otimes dr\otimes1,
\]
is well defined: if $r,r'\in\mathcal J$, then
$d(rr')=r\,dr'+(dr)r'$, and both terms become zero after tensoring on the
left and right with $\mathcal B$.  The quotient map
$q:\mathcal T\to\mathcal B$ induces
\[
\theta:\mathcal B\otimes_{\mathcal T}\Omega^1_{\mathcal T/\mathcal O_X}
\otimes_{\mathcal T}\mathcal B
\longrightarrow\Omega^1_{\mathcal B/\mathcal O_X},
\quad
b_0\otimes dt\otimes b_1\longmapsto b_0\,d(q(t))\,b_1.
\]
The map $\theta$ is surjective because $\Omega^1_{\mathcal B/\mathcal O_X}$
is generated as a $\mathcal B$-bimodule by the elements $db$.  Moreover,
$\theta\circ\partial=0$.

Let $\mathcal C=\operatorname{coker}(\partial)$.  For every
$\mathcal B$-bimodule $\mathcal M$, a morphism
$\mathcal C\to\mathcal M$ is, by the universal property of differentials,
the same as an $\mathcal O_X$-derivation
$\delta:\mathcal T\to\mathcal M$ whose restriction to $\mathcal J$ is
zero.  Such derivations factor uniquely through $\mathcal B$, so
\[
\mathcal Hom_{\mathcal B^e}(\mathcal C,\mathcal M)
\cong\Der_{\mathcal O_X}(\mathcal B,\mathcal M)
\cong\mathcal Hom_{\mathcal B^e}
(\Omega^1_{\mathcal B/\mathcal O_X},\mathcal M)
\]
naturally in $\mathcal M$.  Yoneda identifies $\mathcal C$ with
$\Omega^1_{\mathcal B/\mathcal O_X}$ and proves the right exact sequence.
Applying the left exact functor
$\mathcal Hom_{\mathcal B^e}(-,\mathcal M)$ gives the displayed sequence of
derivations.  Under the explicit identification, the last arrow is
$\delta\mapsto([r]\mapsto\delta(r))$.
\end{proof}

\subsection{The Hochschild degree shift}
There are two distinct infinitesimal questions.  Derivations $B\to B$ are infinitesimal automorphisms and live in Hochschild degree one.  Deformations of the multiplication are Hochschild $2$-cocycles.  Quotienting by infinitesimal changes of basis gives $HH^2$.  The classical first-order deformation complex places these groups in adjacent cohomological degrees:
\[
H^{-1}=\Der_G(B),\qquad
H^0=\frac{\text{$I$-preserving $2$-cocycles}}{\text{coboundaries}}.
\]
Keeping these roles separate prevents the differential calculus from being confused with the moduli tangent space.

\subsection{The PI-restricted tangent complex}

Let $k/F$ be a field extension and let $B=(V_k,\mu)$ be a $k$-point of $\AlgLaws^G_{\mathbf d}(I)$.  We use normalised degree-preserving Hochschild cochains
\[
C^n_{G,\mathrm n}(B,B)=\{\varphi:B^{\otimes n}\to B:\ \varphi\text{ is degree-preserving and }\varphi(a_1,\ldots,a_n)=0\text{ whenever some }a_i=1\}.
\]
We use the following fixed sign convention for the Hochschild differential:
\[
(d_{\mathrm{Hoch}}\psi)(x,y)
=\psi(xy)-\psi(x)y-x\psi(y)
\]
for $\psi\in C^1_{G,\mathrm n}(B,B)$, and
\[
(d_{\mathrm{Hoch}}\varphi)(x,y,z)
=\varphi(xy,z)-\varphi(x,yz)+\varphi(x,y)z-x\varphi(y,z)
\]
for $\varphi\in C^2_{G,\mathrm n}(B,B)$.  This is the negative of another
common convention, but it has the same cocycles, coboundaries, and
cohomology. 
For $f\in I$ and $\varphi\in C^2_{G,\mathrm n}(B,B)$, evaluate $f$ using the first-order multiplication
$
\mu_\varepsilon=\mu+\varepsilon\varphi,
$ with $ \varepsilon^2=0.
$
Since $f$ is an identity of $B$, there is a uniquely determined multilinear map $D_{\mu}f(\varphi)$ such that
$
f_{\mu_\varepsilon}=\varepsilon D_{\mu}f(\varphi).
$

\begin{Definition}\label{def:first-order-PI-deformation-groupoid}
Let $\operatorname{Def}^{(1)}_{G,I}(B)$ be the groupoid whose objects are unital graded multiplications $\mu_\varepsilon$ on $B\otimes_k k[\varepsilon]/(\varepsilon^2)$ reducing to $\mu$ modulo $\varepsilon$ and satisfying $I$, and whose arrows are graded algebra isomorphisms reducing to the identity modulo $\varepsilon$.
\end{Definition}

\begin{Definition}
The space of $I$-preserving graded $2$-cocycles is
\[
Z^2_{G,I,\mathrm n}(B,B)
=
\{\varphi\in C^2_{G,\mathrm n}(B,B):
 d_{\mathrm{Hoch}}\varphi=0,
\ D_\mu f(\varphi)=0\text{ for every }f\in I\}.
\]
\end{Definition}

\begin{Lemma}\label{lem:linearized-generators}
It is enough in the preceding definition to impose $D_\mu f_j(\varphi)=0$ for a fixed $T_G$-generating system $f_1,\ldots,f_r$ of $I$.
\end{Lemma}

\begin{proof}
Fix homogeneous bases of $B$.  Every evaluation of a polynomial identity on
basis vectors gives a vector of polynomial functions on the affine space of
multiplications.  Let $\mathfrak a(f_1,\ldots,f_r)$ be the ideal generated
by all coefficient functions obtained from all graded substitutions of the
$f_j$.  By the definition of $T_G$-generation, every coefficient function
coming from an element $f\in I$ belongs to this ideal.
Let $c$ be such a coefficient function.  Locally in the polynomial ring we
may write
\[
c=\sum_\nu h_\nu c_\nu,
\]
where each $c_\nu$ is a coefficient of a substituted evaluation of one of
the $f_j$.  At the point $\mu$, all $c_\nu(\mu)$ vanish because $B$
satisfies $I$.  Differentiating in the direction $\varphi$ gives
$
(dc)_\mu(\varphi)=
\sum_\nu h_\nu(\mu)(dc_\nu)_\mu(\varphi).
$
The differential of substitution is obtained by evaluating the first
variation $D_\mu f_j(\varphi)$ after the same substitution.  Hence
$D_\mu f_j(\varphi)=0$ for all $j$ implies
$D_\mu f(\varphi)=0$ for every $f\in I$.  The converse follows because
$f_j\in I$.
\end{proof}

\subsection{The deformation groupoid as an action groupoid}
Let $D=k[\varepsilon]/(\varepsilon^2)$. A first-order deformation is a multiplication 
$
\mu_\varphi=\mu+\varepsilon\varphi
$
on $B_D=B\otimes_kD$. A gauge transformation reducing to the identity is
\[
g_\psi=\id+\varepsilon\psi,
\qquad \psi\in C^1_{G,\mathrm n}(B,B).
\]
Transporting $\mu_\varphi$ by $g_\psi$ changes $\varphi$ by a Hochschild coboundary. With the convention
$
\psi\cdot\varphi=\varphi+d_{\mathrm{Hoch}}\psi,
$
this defines an additive action on $Z^2_{G,I,\mathrm n}(B,B)$.

\begin{Lemma}\label{lem:gauge-calculation}
Let $\varphi,\varphi'\in C^2_{G,\mathrm n}(B,B)$ and
$\psi\in C^1_{G,\mathrm n}(B,B)$.  The map
$
g_\psi=\id+\varepsilon\psi
$
is a graded unital algebra isomorphism
$(B_D,\mu+\varepsilon\varphi)\to
(B_D,\mu+\varepsilon\varphi')$ reducing to the identity modulo
$\varepsilon$ if and only if
$
\varphi'=\varphi+d_{\mathrm{Hoch}}\psi.
$
Moreover $g_\chi\circ g_\psi=g_{\psi+\chi}$ and
$g_\psi^{-1}=g_{-\psi}$.
\end{Lemma}
\begin{proof}
The unit and grading conditions are exactly $\psi(1)=0$ and degree
preservation.  Multiplicativity means
\[
g_\psi\bigl(\mu(x,y)+\varepsilon\varphi(x,y)\bigr)
=(\mu+\varepsilon\varphi')
(g_\psi x,g_\psi y).
\]
The coefficient of $\varepsilon$ on the left is
$\varphi(x,y)+\psi(xy)$, while on the right it is
$\varphi'(x,y)+\psi(x)y+x\psi(y)$.  Equality is therefore equivalent to
\[
\varphi'(x,y)=\varphi(x,y)+\psi(xy)-\psi(x)y-x\psi(y),
\]
which is the displayed formula.  The composition and inverse statements
follow by multiplying the endomorphisms modulo $\varepsilon^2$.
\end{proof}

Choose generators $f_1,\dots,f_r$ of $I$. Their first variations define
\[
\Lambda_{I,B}:Z^2_{G,\mathrm n}(B,B)
\longrightarrow
\bigoplus_{j=1}^r\operatorname{Hom}_k(B^{\otimes m_j},B),
\qquad
\varphi\longmapsto(D_\mu f_j(\varphi))_j.
\]
Then $Z^2_{G,I,\mathrm n}(B,B)=\ker\Lambda_{I,B}$.  Since a change of
basis preserves every polynomial identity,
$d_{\mathrm{Hoch}}C^1_{G,\mathrm n}(B,B)\subseteq\ker\Lambda_{I,B}$, and
\[
0\longrightarrow Z^2_{G,I,\mathrm n}(B,B)
\longrightarrow Z^2_{G,\mathrm n}(B,B)
\xrightarrow{\ \Lambda_{I,B}\ }
\bigoplus_{j=1}^r\operatorname{Hom}_k(B^{\otimes m_j},B)
\]
is exact at the first two terms.  Thus the PI-restricted deformation complex
is the subcomplex obtained by imposing the linearised identity equations.

\begin{Theorem}\label{thm:tangent-complex}
Put
$
X_I=\AlgLaws^G_{\mathbf d}(I),\ 
\mathcal X_I=\PI^G_{I,\mathbf d}=[X_I/H_{\mathbf d}],
$
and let $x:\Spec k\to X_I$ correspond to the graded algebra $B=(V_k,\mu)$.  For the induced point $[B]:\Spec k\to\mathcal X_I$, set
\[
\mathbb T_{\mathcal X_I/F,[B]}
:=\operatorname{RHom}_k
\bigl(\mathbb L_{\mathcal X_I/F}\otimes^{\mathbf L}_{\mathcal O_{\mathcal X_I}}k,k\bigr).
\]
There is a canonical quasi-isomorphism
\begin{equation}\label{eq:intrinsic-tangent-complex}
\tau_{\leq0}\mathbb T_{\mathcal X_I/F,[B]}
\simeq
\left[
C^1_{G,\mathrm n}(B,B)
\xrightarrow{d_{\mathrm{Hoch}}}
Z^2_{G,I,\mathrm n}(B,B)
\right],
\end{equation}
where the terms on the right occur in degrees $-1$ and $0$.  We denote the two-term complex on the right by
$\mathbb T^{\mathrm{cl}}_{I,[B]}$.  Its action groupoid is naturally equivalent to
$\operatorname{Def}^{(1)}_{G,I}(B)$.  Consequently,
\[
H^{-1}\mathbb T_{\mathcal X_I/F,[B]}\cong\Der_G(B),
\qquad
H^0\mathbb T_{\mathcal X_I/F,[B]}
\cong
\frac{Z^2_{G,I,\mathrm n}(B,B)}{B^2_{G,\mathrm n}(B,B)},
\]
and the latter vector space embeds naturally in $\HH^2_G(B,B)$.  If $X_I$ is smooth at $x$, then $\mathbb T_{\mathcal X_I/F,[B]}$ has amplitude $[-1,0]$, and the truncation in \eqref{eq:intrinsic-tangent-complex} may be omitted.
\end{Theorem}

\begin{proof}
Let $H=H_{\mathbf d}$ and $\mathfrak h=\Lie(H)$.  The atlas
$p:X_I\to\mathcal X_I$ is an $H$-torsor.  Hence
\[
\mathbb L_{X_I/\mathcal X_I}\otimes^{\mathbf L}k
\simeq\mathfrak h^\vee
\]
in degree $0$.  The transitivity triangle for
$X_I\xrightarrow{p}\mathcal X_I\to\Spec F$, pulled back by $x$ and dualised over $k$, gives
\begin{equation}\label{eq:quotient-tangent-triangle}
\mathfrak h\longrightarrow
\mathbb T_{X_I/F,x}\longrightarrow
\mathbb T_{\mathcal X_I/F,[B]}
\longrightarrow\mathfrak h[1].
\end{equation}
Since $X_I$ is an ordinary finite-type scheme,
$\tau_{\leq0}\mathbb T_{X_I/F,x}=T_xX_I[0]$.  Thus the non-positive truncation of the third term in \eqref{eq:quotient-tangent-triangle} is represented by
\[
[\mathfrak h\longrightarrow T_xX_I]
\]
in degrees $[-1,0]$. 
The Lie algebra $\mathfrak h$ consists of degree-preserving endomorphisms
$\psi$ of $B$ satisfying $\psi(1)=0$; hence
\begin{equation}\label{eq:lie-C1}
\mathfrak h=C^1_{G,\mathrm n}(B,B).
\end{equation}
A tangent vector to $X_I$ at $x$ is a law
$\mu+\varepsilon\varphi$ over $k[\varepsilon]/(\varepsilon^2)$.  The two unit equations give
$\varphi(1,-)=\varphi(-,1)=0$.  The coefficient of $\varepsilon$ in the associator is
$d_{\mathrm{Hoch}}\varphi$, and the coefficient of $\varepsilon$ in the evaluation of $f\in I$ is
$D_\mu f(\varphi)$.  Therefore
\begin{equation}\label{eq:tangent-law-Z2}
T_xX_I=Z^2_{G,I,\mathrm n}(B,B).
\end{equation}
The action of $H$ on laws is
$g\cdot\mu=g\circ\mu\circ(g^{-1}\otimes g^{-1})$.  Substituting
$g=\id+\varepsilon\psi$ shows that its differential at the identity is
\begin{equation}\label{eq:infinitesimal-orbit-Hochschild}
\psi\longmapsto
\bigl((x,y)\mapsto\psi(xy)-\psi(x)y-x\psi(y)\bigr)
=d_{\mathrm{Hoch}}\psi.
\end{equation}
Equations \eqref{eq:lie-C1}--\eqref{eq:infinitesimal-orbit-Hochschild} identify the complex obtained from \eqref{eq:quotient-tangent-triangle} with the right-hand side of \eqref{eq:intrinsic-tangent-complex}.

By Lemma~\ref{lem:gauge-calculation}, an arrow between
$\mu+\varepsilon\varphi$ and $\mu+\varepsilon\varphi'$ reducing to the identity is uniquely $\id+\varepsilon\psi$, and it exists precisely when
$\varphi'=\varphi+d_{\mathrm{Hoch}}\psi$.  This proves the equivalence with the action groupoid.  Its automorphisms are the $1$-cocycles, namely $\Der_G(B)$, and its isomorphism classes form
$Z^2_{G,I,\mathrm n}/B^2_{G,\mathrm n}$.  Since every change-of-basis orbit is contained in $X_I$, one has
$B^2_{G,\mathrm n}\subseteq Z^2_{G,I,\mathrm n}\subseteq Z^2_{G,\mathrm n}$, which gives the asserted injection into $\HH^2_G(B,B)$.

If $X_I$ is smooth at $x$, then
$\mathbb L_{X_I/F}\otimes k$ is concentrated in degree $0$.  The triangle dual to the transitivity triangle therefore shows that
$\mathbb T_{\mathcal X_I/F,[B]}$ is concentrated in degrees $[-1,0]$, proving the last assertion.
\end{proof}

\subsection{Additional identities and the normal complex}

Let $I\subseteq J$ be finitely generated graded $T$-ideals and let
$B=(V_k,\mu)$ satisfy $J$.  Write
$
X_J\xhookrightarrow{i}X_I,
 $ and $
\mathcal X_J\xhookrightarrow{\bar i}\mathcal X_I
$
for the closed immersions of law schemes and quotient stacks.  Let
$\mathcal K_{J/I}\subseteq\mathcal O_{X_I}$ be the ideal sheaf of $X_J$ and put
\[
\mathcal C_{J/I,B}
=(\mathcal K_{J/I}/\mathcal K_{J/I}^2)\otimes_{\mathcal O_{X_J}}k,
\qquad
N_{J/I,B}=\operatorname{Hom}_k(\mathcal C_{J/I,B},k).
\]

\begin{Theorem}\label{thm:relative-PI-tangent}
Under the identifications of Theorem~\ref{thm:tangent-complex}, the morphism induced by $\bar i$ on non-positive tangent complexes is represented by the commutative diagram
\begin{equation}\label{eq:relative-tangent-square}
\begin{tikzcd}[column sep=large]
C^1_{G,\mathrm n}(B,B) \arrow[r,"d_{\mathrm{Hoch}}"] \arrow[d,equal]
& Z^2_{G,J,\mathrm n}(B,B) \arrow[d,hook]\\
C^1_{G,\mathrm n}(B,B) \arrow[r,"d_{\mathrm{Hoch}}"']
& Z^2_{G,I,\mathrm n}(B,B).
\end{tikzcd}
\end{equation}
The cone of \eqref{eq:relative-tangent-square} is canonically quasi-isomorphic to
\begin{equation}\label{eq:first-normal-quotient}
\frac{Z^2_{G,I,\mathrm n}(B,B)}{Z^2_{G,J,\mathrm n}(B,B)}[0].
\end{equation}
Moreover, the conormal sequence induces an exact sequence
\begin{equation}\label{eq:normal-linearisation-sequence}
0\longrightarrow Z^2_{G,J,\mathrm n}(B,B)
\longrightarrow Z^2_{G,I,\mathrm n}(B,B)
\xrightarrow{\lambda_{J/I,B}}N_{J/I,B},
\end{equation}
where, for a local coefficient equation $Q\in\mathcal K_{J/I}$,
\begin{equation}\label{eq:normal-map-equations}
\lambda_{J/I,B}(\varphi)([Q])=(dQ)_B(\varphi).
\end{equation}
If $q_1,\ldots,q_s$ are additional graded identities generating $J$ modulo $I$, then \eqref{eq:normal-map-equations} is computed by the coefficient functions of
 $
D_\mu q_1(\varphi),\ldots,D_\mu q_s(\varphi).
$
In particular,
\begin{equation}\label{eq:normal-image}
Z^2_{G,I,\mathrm n}(B,B)/Z^2_{G,J,\mathrm n}(B,B)
\cong\operatorname{im}\lambda_{J/I,B}.
\end{equation}
\end{Theorem}

\begin{proof}
Both quotient stacks use the same change-of-basis group $H_{\mathbf d}$.  Hence the map of quotient tangent complexes is the identity on
$\Lie(H_{\mathbf d})=C^1_{G,\mathrm n}(B,B)$ and the differential of
$i:X_J\hookrightarrow X_I$ in degree $0$.  By \eqref{eq:tangent-law-Z2}, this differential is the inclusion
$Z^2_{G,J,\mathrm n}\hookrightarrow Z^2_{G,I,\mathrm n}$, which proves \eqref{eq:relative-tangent-square}.  Cancelling the identical degree $-1$ terms gives \eqref{eq:first-normal-quotient}.
The conormal sequence of the closed immersion $i$ is right exact:
\[
\mathcal K_{J/I}/\mathcal K_{J/I}^2
\longrightarrow i^*\Omega^1_{X_I/F}
\longrightarrow\Omega^1_{X_J/F}\longrightarrow0.
\]
Tensoring with $k$ and applying $\operatorname{Hom}_k(-,k)$ yields an exact sequence
\[
0\longrightarrow T_BX_J\longrightarrow T_BX_I
\xrightarrow{\lambda_{J/I,B}}N_{J/I,B}.
\]
Substitution of \eqref{eq:tangent-law-Z2} for $I$ and $J$ gives \eqref{eq:normal-linearisation-sequence}.  By the definition of the dual conormal map, a tangent vector $\varphi$ acts on the class of a local equation $Q$ by directional differentiation, proving \eqref{eq:normal-map-equations}.  The coefficient equations obtained by evaluating the $q_\nu$ generate $\mathcal K_{J/I}$; their directional derivatives are exactly the coefficients of $D_\mu q_\nu(\varphi)$.  Finally, exactness of \eqref{eq:normal-linearisation-sequence} gives \eqref{eq:normal-image}.
\end{proof}

\begin{Corollary}\label{cor:PI-normal-codimension}
If $X_I$ and $X_J$ are smooth at $B$, then
$\mathcal X_J\hookrightarrow\mathcal X_I$ is a regular immersion at $[B]$ and
\begin{equation}\label{eq:PI-codimension}
\operatorname{codim}_{[B]}(\mathcal X_J,\mathcal X_I)
=
\dim_k Z^2_{G,I,\mathrm n}(B,B)
-
\dim_k Z^2_{G,J,\mathrm n}(B,B)
=
\operatorname{rank}\lambda_{J/I,B}.
\end{equation}
More generally, suppose that $X_I$ is smooth at $B$ and that, near $B$, the ideal of $X_J$ in $X_I$ is generated by $c$ coefficient equations whose differentials are linearly independent at $B$.  Then $X_J$ is smooth at $B$, and both the scheme and stack immersions have codimension $c$ there.
\end{Corollary}

\begin{proof}
A closed immersion between schemes smooth over $k$ is regular.  Its dual conormal sequence is therefore short exact, so
$\lambda_{J/I,B}$ is surjective and its rank equals the codimension.  The two quotient stacks are presented by the same smooth group $H_{\mathbf d}$; regularity and codimension may be checked after the common smooth atlas, giving \eqref{eq:PI-codimension}.  The final assertion is the Jacobian criterion applied inside the smooth scheme $X_I$; the same smooth-local argument transfers it to the quotient stacks.
\end{proof}

\begin{Definition}
The algebra $B$ is $I$-rigid if $H^0(\mathbb T^{\mathrm{cl}}_{I,[B]})=0$.
\end{Definition}

\begin{Corollary}\label{cor:rigidity-orbit}
If $B$ is $I$-rigid and $\AlgLaws^G_{\mathbf d}(I)$ is smooth at $B$, then the $H_{\mathbf d}$-orbit of $B$ is open in a neighbourhood of $B$.
\end{Corollary}

\begin{proof}
The orbit is smooth because $H_{\mathbf d}$ and the stabiliser of $B$ are
smooth in characteristic zero.  Rigidity identifies its tangent space at
$B$ with $T_B\AlgLaws^G_{\mathbf d}(I)$.  Since the law scheme is smooth at
$B$, the locally closed immersion of the orbit is \'{e}tale at $B$; after
shrinking around $B$, it is therefore an open immersion.
\end{proof}

\subsection{Explicit linearisations}
For the associator $a(x,y,z)=(xy)z-x(yz)$,
\[
D_\mu a(\varphi)(x,y,z)=
\varphi(xy,z)-\varphi(x,yz)+\varphi(x,y)z-x\varphi(y,z)
=d_{\mathrm{Hoch}}\varphi(x,y,z).
\]
For the commutator identity,
$
D_\mu[x,y](\varphi)=\varphi(x,y)-\varphi(y,x).
$

\begin{Proposition}\label{prop:associative-harrison}
Assume that $G$ is trivial.
\begin{enumerate}
\item If no identities are imposed beyond associativity and the unit, then
$
H^0\mathbb T_{\PI_{0,d}/F,[B]}\cong\HH^2(B,B).
$
\item If $B$ is commutative and $J$ is the $T$-ideal generated by $[x_1,x_2]$, then
$
H^0\mathbb T_{\PI_{J,d}/F,[B]}\cong\operatorname{Harr}^2(B,B).
$
\end{enumerate}
\end{Proposition}

\begin{proof}
For the first assertion,
$Z^2_{0,\mathrm n}(B,B)=Z^2_{\mathrm n}(B,B)$, and Theorem~\ref{thm:tangent-complex} gives
$H^0=Z^2_{\mathrm n}/B^2_{\mathrm n}=\HH^2(B,B)$, the classical tangent space to associative deformations \cite{Gerstenhaber}. 
For the second assertion, the displayed linearisation of the commutator shows that
$Z^2_{J,\mathrm n}(B,B)$ is the space of symmetric normalised Hochschild $2$-cocycles.  Since $B$ is commutative, every normalised Hochschild coboundary is symmetric.  In characteristic zero, the degree-$2$ Harrison complex is precisely the symmetric summand of the normalised Hochschild complex.  Therefore
\[
Z^2_{J,\mathrm n}(B,B)/B^2_{\mathrm n}(B,B)
=\operatorname{Harr}^2(B,B),
\]
as asserted \cite{Harrison}.
\end{proof}

In the supercommutative case, the commutator linearisation is replaced by its Koszul-signed analogue.

\subsection{Second-order PI obstructions}

The preceding tangent complex records which first-order directions satisfy the
linearised equations.  The following result determines when such a direction extends one
order further.  Put
\[
D_2=k[t]/(t^2),\qquad D_3=k[t]/(t^3),
\]
and let $\varphi\in Z^2_{G,I,\mathrm n}(B,B)$.  A second-order lift of
$\mu+t\varphi$ is a multiplication
\[
\mu_{\varphi,\psi}=\mu+t\varphi+t^2\psi
\quad\text{on }B\otimes_kD_3,
\qquad
\psi\in C^2_{G,\mathrm n}(B,B),
\]
which is associative, unital, graded, satisfies $I$, and reduces to
$\mu+t\varphi$ over $D_2$. 
Define the quadratic associativity term
\[
Q_{\mathrm{ass}}(\varphi)(x,y,z)
=
\varphi(\varphi(x,y),z)-\varphi(x,\varphi(y,z)).
\]
For a normalised $2$-cochain $\psi$, write
\[
\partial_\mu\psi(x,y,z)
=
\psi(xy,z)-\psi(x,yz)+\psi(x,y)z-x\psi(y,z).
\]
Thus $\partial_\mu=d_{\mathrm{Hoch}}$ on normalised $2$-cochains for
the convention fixed above.

Fix finite multilinear graded $T_G$-generators
$f_1,\ldots,f_r$ of $I$.  Using the fixed bracketings chosen in the
defining equations of the law scheme, evaluation of $f_j$ is polynomial in
the multiplication law, so its second differential at $\mu$ is defined.  Set
\[
Q_j(\varphi)=\frac12D^2_\mu f_j(\varphi,\varphi).
\]
Equivalently, because $f_j(\mu)=0$ and $D_\mu f_j(\varphi)=0$,
$Q_j(\varphi)$ is the coefficient of $t^2$ in the evaluation of $f_j$
using $\mu+t\varphi$.
Consider the linear operator
\[
\mathcal L_{I,B}:C^2_{G,\mathrm n}(B,B)
\longrightarrow
C^3_{G,\mathrm n}(B,B)
\oplus
\bigoplus_{j=1}^r\operatorname{Hom}_k(B^{\otimes m_j},B)
\]
given by
$
\mathcal L_{I,B}(\psi)
=
\bigl(\partial_\mu\psi,
D_\mu f_1(\psi),\ldots,D_\mu f_r(\psi)\bigr),
$
where $m_j$ is the number of variables of $f_j$.  The second-order obstruction associated with this presentation is
\[
\operatorname{ob}^{(2)}_{I,B}(\varphi)
=
\left[
Q_{\mathrm{ass}}(\varphi),
Q_1(\varphi),\ldots,Q_r(\varphi)
\right]
\in\operatorname{coker}(\mathcal L_{I,B}).
\]

\begin{Theorem}\label{thm:second-order-PI-obstruction}
For $\varphi\in Z^2_{G,I,\mathrm n}(B,B)$, the following conditions are
equivalent.
\begin{enumerate}
\item The first-order deformation $\mu+t\varphi$ over $D_2$ lifts to an
$I$-preserving unital graded multiplication over $D_3$.
\item There exists $\psi\in C^2_{G,\mathrm n}(B,B)$ satisfying
$
\partial_\mu\psi=-Q_{\mathrm{ass}}(\varphi),
 $ and $
D_\mu f_j(\psi)=-Q_j(\varphi)
\quad(1\leq j\leq r).
$
\item $\operatorname{ob}^{(2)}_{I,B}(\varphi)=0$.
\end{enumerate}
The displayed cokernel depends on the chosen generators.  Its zero condition
does not: it is equivalent to liftability, hence is independent of the
generating system and invariant under isomorphisms of first-order
deformations.
\end{Theorem}

\begin{proof}
Let
$
\mu_t=\mu+t\varphi+t^2\psi.
$
We first calculate the associator.  Expanding bilinearly and collecting
powers of $t$ gives
\[
\mu_t(\mu_t(x,y),z)-\mu_t(x,\mu_t(y,z))
=t\,\partial_\mu\varphi(x,y,z)
+t^2\!\left(\partial_\mu\psi(x,y,z)
+\varphi(\varphi(x,y),z)-\varphi(x,\varphi(y,z))\right)
\pmod{t^3}.
\]
The constant term is zero because $\mu$ is associative, and the coefficient
of $t$ is zero because
$\varphi\in Z^2_{G,I,\mathrm n}(B,B)$.  Hence $\mu_t$ is associative modulo
$t^3$ if and only if
$
\partial_\mu\psi=-Q_{\mathrm{ass}}(\varphi).
$
Normalisation of $\varphi$ and $\psi$ is equivalent to the two unit
equations modulo $t^3$, and the grading is preserved by construction.
For each generator $f_j$, evaluation is a polynomial map from the affine
space of multiplication laws to
$\operatorname{Hom}_k(B^{\otimes m_j},B)$.  Since $F$ has characteristic
zero, Taylor expansion along the line
$\mu+t\varphi+t^2\psi$ is valid and yields
\[
f_j(\mu_t)=f_j(\mu)+tD_\mu f_j(\varphi)
+t^2\!\left(D_\mu f_j(\psi)+\frac12D^2_\mu f_j(\varphi,\varphi)\right)
\pmod{t^3}.
\]
The first two terms vanish because $B$ satisfies $I$ and $\varphi$ is an
$I$-preserving first-order cocycle.  Therefore $f_j(\mu_t)=0$ modulo $t^3$
if and only if
\[
D_\mu f_j(\psi)=-Q_j(\varphi).
\]
Because the $f_j$ generate $I$ as a graded $T$-ideal, vanishing of their
evaluations over the commutative base $D_3$ is equivalent to vanishing of
every element of $I$.  This proves the equivalence of (1) and (2).
Condition (2) says exactly
\[
\mathcal L_{I,B}(\psi)
=-\bigl(Q_{\mathrm{ass}}(\varphi),Q_1(\varphi),\ldots,Q_r(\varphi)\bigr).
\]
Since the image of a linear map is closed under multiplication by $-1$,
this is equivalent to the vanishing of the displayed cokernel class, proving
(2)$\Longleftrightarrow$(3).

It remains to justify the two intrinsicity assertions.  If another finite generating
system is chosen, the corresponding equations define the same functor of
$D_3$-valued algebra laws: a multiplication satisfies either system if and
only if it satisfies the $T_G$-ideal $I$.  The existence of $\psi$ solving
the lifting problem is therefore independent of the presentation.  The
actual representatives in the two cokernels need not be canonically equal;
what is intrinsic is their vanishing, equivalently the existence of a lift.

Finally suppose
$g_2:(B_{D_2},\mu+t\varphi)\xrightarrow{\sim}
(B_{D_2},\mu+t\varphi')$
is an isomorphism reducing to the identity.  It has the form
$g_2=\id+t\alpha$.  The same formula defines an invertible graded
$D_3$-linear map $\widetilde g=\id+t\alpha$ with inverse
$\id-t\alpha+t^2\alpha^2$.  If $\widetilde\mu$ is a second-order lift of the
first deformation, transport it by $\widetilde g$:
\[
\widetilde\mu'(x,y)=
\widetilde g\bigl(
\widetilde\mu(\widetilde g^{-1}x,\widetilde g^{-1}y)
\bigr).
\]
Transport preserves associativity, the unit, the grading, and every
polynomial identity, and reduction modulo $t^2$ is the target deformation.
Thus liftability is invariant under isomorphism.  Applying the same argument
to $g_2^{-1}$ proves the converse.
\end{proof}

\begin{Corollary}\label{cor:smooth-formal-liftability}
If $B$ is a smooth point of $\AlgLaws^G_{\mathbf d}(I)$, then every
$I$-preserving first-order deformation of $B$ lifts to a compatible unital
graded $I$-preserving multiplication over $k[t]/(t^{N+1})$ for every
$N\geq2$.  In particular,
$\operatorname{ob}^{(2)}_{I,B}(\varphi)=0$ for every
$\varphi\in Z^2_{G,I,\mathrm n}(B,B)$.
\end{Corollary}

\begin{proof}
A first-order deformation is a morphism
$\operatorname{Spec}D_2\to\AlgLaws^G_{\mathbf d}(I)$ whose closed point is
$B$.  Smoothness at $B$ gives the infinitesimal lifting property after
restricting to a neighbourhood of $B$.  Apply it successively to the
square-zero surjections
\[
k[t]/(t^{n+2})\twoheadrightarrow k[t]/(t^{n+1}),
\qquad n\geq1.
\]
This constructs compatible lifts to every finite order.  The final assertion
follows from Theorem~\ref{thm:second-order-PI-obstruction}.
\end{proof}

\begin{Corollary}\label{cor:linear-second-order-criterion}
If $\mathcal L_{I,B}$ is surjective, then every element of
$Z^2_{G,I,\mathrm n}(B,B)$ lifts to second order.
\end{Corollary}

\begin{proof}
If $\mathcal L_{I,B}$ is surjective, every obstruction vector lies in its image; the assertion follows from Theorem~\ref{thm:second-order-PI-obstruction}.
\end{proof}

\subsection{Relative K-theory of PI-preserving first-order deformations}\label{subsec:PI-infinitesimal-K}

Let $k/F$ be a field extension and let $B$ be a finite-dimensional unital graded $k$-algebra satisfying $I$.  For $\varphi\in Z^2_{G,I,\mathrm n}(B,B)$, set
$
B_\varphi=B\oplus\varepsilon B, $with $ \varepsilon^2=0,
$
with multiplication
\begin{equation}\label{eq:PI-square-zero-product}
(a+\varepsilon x)\star_\varphi(b+\varepsilon y)
=ab+\varepsilon\bigl(ay+xb+\varphi(a,b)\bigr).
\end{equation}
The projection $q_\varphi:B_\varphi\to B$ is an augmented algebra homomorphism with square-zero kernel $\varepsilon B$.  Write
\[
\mathbf K(B_\varphi,B)
:=\operatorname{fib}\bigl(\mathbf K(B_\varphi)\longrightarrow\mathbf K(B)\bigr)
\]
for the relative connective K-theory spectrum.  We write $\mathbf{HN}(B_\varphi,B)$ for the Eilenberg--Mac Lane spectrum of the relative negative cyclic complex, indexed so that
$\pi_m\mathbf{HN}(B_\varphi,B)=HN_m(B_\varphi,B)$.

\begin{Theorem}\label{thm:PI-relative-K-functor}
The assignment
\[
\varphi\longmapsto \mathbf K(B_\varphi,B)
\]
extends canonically to a functor
\[
\mathbf K^{\mathrm{rel}}_{G,I,B}:
\operatorname{Def}^{(1)}_{G,I}(B)\longrightarrow\mathbf{Sp}^{\simeq}.
\]
Hence isomorphic PI-preserving first-order deformations determine equivalent relative K-theory spectra.  Under the action-groupoid description of Theorem~\ref{thm:tangent-complex}, the arrow represented by $\psi\in C^1_{G,\mathrm n}(B,B)$ is induced by the augmented algebra isomorphism
\[
g_\psi=\id+\varepsilon\psi:
(B_\varphi,q_\varphi)\xrightarrow{\sim}
(B_{\varphi+d_{\mathrm{Hoch}}\psi},q_{\varphi+d_{\mathrm{Hoch}}\psi}).
\]
\end{Theorem}

\begin{proof}
The coefficient of $\varepsilon$ in the associator of \eqref{eq:PI-square-zero-product} is $d_{\mathrm{Hoch}}\varphi$, while the coefficient of $\varepsilon$ in the evaluation of $f\in I$ is $D_\mu f(\varphi)$.  Thus $B_\varphi$ is associative, unital, graded and satisfies $I$ precisely because $\varphi\in Z^2_{G,I,\mathrm n}(B,B)$.  Lemma~\ref{lem:gauge-calculation} shows that $g_\psi$ is an algebra isomorphism with the stated target.  Since $g_\psi$ is the identity modulo $\varepsilon$, one has
$q_{\varphi+d_{\mathrm{Hoch}}\psi}\circ g_\psi=q_\varphi$.
Consequently it induces a morphism between the homotopy fibres defining relative K-theory.  Moreover,
$g_{\psi'}g_\psi=g_{\psi+\psi'}$ because $\varepsilon^2=0$, and $g_0=\id$.  Functoriality follows.
\end{proof}

\begin{Theorem}\label{thm:Goodwillie-PI-relative}
For every $\varphi\in Z^2_{G,I,\mathrm n}(B,B)$ and every $m\geq1$, the relative Chern character induces a natural isomorphism
\begin{equation}\label{eq:Goodwillie-PI-comparison}
\operatorname{ch}^{-}_m:
K_m(B_\varphi,B)\otimes_{\mathbb Z}\mathbb Q
\xrightarrow{\ \sim\ }
HN_m(B_\varphi,B).
\end{equation}
Equivalently, with the indexing convention above, there is a natural equivalence
\[
\tau_{\geq1}\bigl(\mathbf K(B_\varphi,B)\otimes\mathbb Q\bigr)
\simeq
\tau_{\geq1}\mathbf{HN}(B_\varphi,B).
\]
These equivalences are natural on $\operatorname{Def}^{(1)}_{G,I}(B)$.
\end{Theorem}

\begin{proof}
The field $k$ is a $\mathbb Q$-algebra and the kernel $\varepsilon B$ of $q_\varphi$ satisfies $(\varepsilon B)^2=0$.  Goodwillie's nilpotent comparison theorem therefore identifies relative algebraic K-theory rationally with relative negative cyclic homology.  The relative Chern character gives the isomorphism in \eqref{eq:Goodwillie-PI-comparison}; the identification with Goodwillie's rational homotopy character is proved by Corti\~nas and Weibel \cite{GoodwillieRelative,CortinasWeibelRelative}.  Both constructions are natural for morphisms of augmented algebras.  Applying this naturality to the maps $g_\psi$ of Theorem~\ref{thm:PI-relative-K-functor} proves the final assertion.
\end{proof}

\begin{Example}\label{ex:dual-numbers-relative-K}
Assume that $G$ is trivial, let $I$ be generated by $[x_1,x_2]$, and take
$B=k[u]/(u^2)$.  Define the normalised symmetric cochain $\varphi$ by
$\varphi(u,u)=1$ and by zero on all pairs containing $1$.  The only nontrivial value of the Hochschild cocycle equation is
\[
(d_{\mathrm{Hoch}}\varphi)(u,u,u)
=\varphi(u^2,u)-\varphi(u,u^2)+\varphi(u,u)u-u\varphi(u,u)=0,
\]
so $\varphi\in Z^2_{I,\mathrm n}(B,B)$.  In $B_\varphi$ one has $u^2=\varepsilon$; hence
$
B_\varphi\cong k[T]/(T^4) $ and $
q_\varphi:k[T]/(T^4)\longrightarrow k[T]/(T^2).
$
Let $J=(T^2)/(T^4)$.  The algebra $k[T]/(T^4)$ is commutative and local, and $J$ lies in its maximal ideal.  The determinant description of relative $K_1$ for a commutative local ring therefore gives
\[
K_1(B_\varphi,B)
\cong\ker\bigl((k[T]/(T^4))^\times\to(k[T]/(T^2))^\times\bigr)
=1+J.
\]
Since $J^2=0$, the map $1+x\mapsto x$ is an isomorphism from $1+J$ to the additive group of $J$.  Moreover, every normalised Hochschild coboundary takes $(u,u)$ to an element of $ku$, so the class of $\varphi$, whose value there is $1$, is nonzero.  Therefore
\[
K_1(B_\varphi,B)\cong J\cong kT^2\oplus kT^3.
\]
Thus a nonzero Harrison direction produces a nontrivial relative K-group, and Theorem~\ref{thm:Goodwillie-PI-relative} identifies its rational Chern character with the corresponding relative negative cyclic class.
\end{Example}

\subsection{Localisation along additional identities}\label{subsec:PI-K-localisation}

Let $I\subseteq J$ be finitely generated graded $T$-ideals.  Put
\[
X_I=\AlgLaws^G_{\mathbf d}(I),\qquad X_J=\AlgLaws^G_{\mathbf d}(J),
\]
and let $X_I^{\Az}\subseteq X_I$ be the inverse image of the Azumaya substack under the atlas $X_I\to\PI^G_{I,\mathbf d}$.  Set
\[
Z=X_J\times_{X_I}X_I^{\Az},\qquad U=X_I^{\Az}\setminus Z,
\]
and denote by $\mathcal A_I$ the tautological Azumaya algebra on $X_I^{\Az}$.  For a locally ringed space $Y$ equipped with an Azumaya algebra $\mathcal A$, let $\Perf(Y;\mathcal A)$ be the small idempotent-complete stable category of perfect left $\mathcal A$-modules, and let $\Perf_Z(Y;\mathcal A)$ be its full stable subcategory of objects whose underlying $\mathcal O_Y$-complex is supported on $Z$.

\begin{Theorem}\label{thm:PI-identity-localisation}
Restriction to $U$ gives an exact sequence of small idempotent-complete stable categories
\begin{equation}\label{eq:PI-perfect-localisation}
\Perf_Z(X_I^{\Az};\mathcal A_I)
\longrightarrow
\Perf(X_I^{\Az};\mathcal A_I)
\longrightarrow
\Perf(U;\mathcal A_I|_U),
\end{equation}
where the third term is the idempotent completion of the Verdier quotient of the second by the first.  Consequently, nonconnective algebraic K-theory gives a fibre sequence
\begin{equation}\label{eq:PI-K-localisation}
\mathbf K^{\mathrm{nc}}_Z(X_I^{\Az};\mathcal A_I)
\longrightarrow
\mathbf K^{\mathrm{nc}}(X_I^{\Az};\mathcal A_I)
\longrightarrow
\mathbf K^{\mathrm{nc}}(U;\mathcal A_I|_U).
\end{equation}
Let
$
\mathcal X=[X_I^{\Az}/H_{\mathbf d}], 
\mathcal Z=[Z/H_{\mathbf d}], $ and $
\mathcal U=[U/H_{\mathbf d}].
$
If $\mathcal X$ is a perfect algebraic stack, in the sense that its derived category of quasi-coherent sheaves is compactly generated and its compact objects are perfect, then the analogous exact sequence and K-theory fibre sequence hold for the tautological Azumaya algebra on
$\mathcal Z\hookrightarrow\mathcal X\hookleftarrow\mathcal U$.
\end{Theorem}

\begin{proof}
The closed immersion $X_J\hookrightarrow X_I$ is defined by the coefficient equations obtained from the additional identities in $J$; hence $Z\hookrightarrow X_I^{\Az}$ is closed and $U$ is its open complement.  The algebra $\mathcal A_I$ is Azumaya, so it is \'{e}tale-locally a matrix algebra.  On every splitting open, perfect $\mathcal A_I$-modules are Morita equivalent to ordinary perfect complexes, and support of an $\mathcal A_I$-module is the support of its underlying $\mathcal O$-complex.  The Thomason--Trobaugh localisation theorem therefore identifies the kernel of restriction to $U$ with $\Perf_Z(X_I^{\Az};\mathcal A_I)$ and identifies the idempotent completion of the Verdier quotient with $\Perf(U;\mathcal A_I|_U)$ \cite{ThomasonTrobaugh}.  These identifications descend from an \'{e}tale splitting cover because both perfect modules and support satisfy descent.  This proves \eqref{eq:PI-perfect-localisation}.

Nonconnective K-theory is a localising invariant of small idempotent-complete stable categories; applying it to \eqref{eq:PI-perfect-localisation} gives \eqref{eq:PI-K-localisation} \cite{BlumbergGepnerTabuada}.  Under the perfectness hypothesis on $\mathcal X$, the same localisation theorem applies to perfect complexes on the quotient stack and, by the same \'{e}tale Morita argument, to perfect modules over its tautological Azumaya algebra.  This proves the stack statement.
\end{proof}

\subsection{Orbit geometry and rigidity}
\begin{Proposition}\label{prop:orbit-normal-space}
Let $\mathcal O_B\subseteq\AlgLaws^G_{\mathbf d}(I)$ be the change-of-basis orbit of $B$.  Then
$
T_B\mathcal O_B=B^2_{G,\mathrm n}(B,B),
$ $
T_B\AlgLaws^G_{\mathbf d}(I)=Z^2_{G,I,\mathrm n}(B,B),
$
and hence
$
T_B\AlgLaws^G_{\mathbf d}(I)/T_B\mathcal O_B
\cong H^0(\mathbb T^{\mathrm{cl}}_{I,[B]}).
$
Moreover,
\[
\operatorname{Lie}\Aut_G(B)=\Der_G(B).
\]
\end{Proposition}
\begin{proof}
Let $\psi\in C^1_{G,\mathrm n}(B,B)$ and consider the curve of basis changes
$g_t=\id+t\psi$ over the dual numbers.  The gauge calculation of
Lemma~\ref{lem:gauge-calculation} shows that the derivative of the orbit map
at the identity is $\psi\mapsto d_{\mathrm{Hoch}}\psi$; hence its image is
$B^2_{G,\mathrm n}(B,B)$.  Linearising the equations defining the law
scheme gives the associativity equation $d_{\mathrm{Hoch}}\varphi=0$ and the
identity equations $D_\mu f(\varphi)=0$, so its tangent space is
$Z^2_{G,I,\mathrm n}(B,B)$.  The quotient is therefore the degree-zero
cohomology of the complex in Theorem~\ref{thm:tangent-complex}.  The kernel
of the infinitesimal orbit map consists of the $1$-cochains satisfying
$d_{\mathrm{Hoch}}\psi=0$, namely the graded derivations.
\end{proof}

\begin{Example}
Over a field of characteristic zero, a finite-dimensional separable algebra has vanishing positive Hochschild cohomology with coefficients in itself.  In particular, $HH^2(M_n(F),M_n(F))=0$, so the matrix point is rigid in every PI-substack containing it.  Its stabiliser is $\operatorname{PGL}_n$ in the ungraded case.
\end{Example}
\section{Azumaya and form loci}\label{sec:azumaya}

Let $\mathcal B$ be a finite locally free algebra over a scheme $S$.  There is a canonical morphism of vector bundles of the same rank
\[
\Phi_{\mathcal B}:\mathcal B\otimes_{\mathcal O_S}\mathcal B^{\op}
\longrightarrow
\cEnd_{\mathcal O_S}(\mathcal B),
\qquad
b\otimes c\longmapsto(x\mapsto bxc).
\]

\begin{Theorem}\label{thm:Azumaya-open}
The locus
\[
S^{\Az}(\mathcal B)=\{s\in S:(\Phi_{\mathcal B})_s\text{ is an isomorphism}\}
\]
is open and is stable under arbitrary base change.  Consequently, $\PI^G_{I,\mathbf d}$ contains an open substack
\[
\PI^{G,\Az}_{I,\mathbf d}\subseteq\PI^G_{I,\mathbf d}
\]
classifying families whose underlying ungraded algebras are Azumaya over the base.
\end{Theorem}

\begin{proof}
Both sides of $\Phi_{\mathcal B}$ are locally free of rank $(\rank\mathcal B)^2$.  The isomorphism locus of a morphism between vector bundles of equal rank is the complement of the vanishing locus of its determinant, hence is open.  Formation of $\Phi_{\mathcal B}$ commutes with base change.  Applying this construction to the universal family on the quotient stack gives an invariant open substack.
\end{proof}

For an Azumaya family, the quotient of all $F$-linear derivations by inner derivations is identified, \'{e}tale-locally, with the derivations of the centre.

\begin{Theorem}\label{thm:Azumaya-derivations}
Let $S$ be an $F$-scheme and let $\mathcal B$ be an Azumaya $\mathcal O_S$-algebra.  Then, on the small \'{e}tale site of $S$, there are canonical exact sequences of sheaves of Lie algebras
\[
0\longrightarrow\mathcal O_S
\longrightarrow\mathcal B
\xrightarrow{\operatorname{ad}}
\Der_{\mathcal O_S}(\mathcal B)
\longrightarrow0
\]
and
\[
0\longrightarrow\mathcal B/\mathcal O_S
\longrightarrow\Der_F(\mathcal B)
\xrightarrow{\mathrm{res}}
\Der_F(\mathcal O_S)
\longrightarrow0.
\]
\end{Theorem}

\begin{proof}
All assertions are local for the small \'{e}tale topology, so choose an \'{e}tale
map $U\to S$ and an algebra isomorphism
$\mathcal B|_U\cong M_n(\mathcal O_U)$.  We first treat a commutative
$F$-algebra $R$ and $A=M_n(R)$.
Let $\delta\in\Der_R(A)$.  With $e_{ij}$ denoting the matrix units, a direct
calculation from
$e_{ij}e_{kl}=\delta_{jk}e_{il}$ shows that
\[
a=\sum_{i=1}^n\delta(e_{i1})e_{1i}
\]
satisfies $\delta(e_{ij})=[a,e_{ij}]$ for all $i,j$.  Since the matrix units
generate $A$ as an $R$-module, $\delta=\operatorname{ad}(a)$.  The kernel
of $a\mapsto\operatorname{ad}(a)$ consists of matrices commuting with every
$e_{ij}$, namely the scalar matrices $R\cdot I_n$.  Hence
\[
0\longrightarrow R\longrightarrow M_n(R)
\xrightarrow{\operatorname{ad}}\Der_R(M_n(R))
\longrightarrow0
\]
is exact.  These local sequences are preserved by change of \'{e}tale
trivialisation because inner derivations are intrinsic, and therefore they
descend to the first exact sequence of sheaves.
Now, let $D\in\Der_F(\mathcal B)$.  If $z$ is central and $b$ is any local
section, then
\[
0=D([z,b])=[D(z),b]+[z,D(b)]=[D(z),b],
\]
so $D(z)$ is central.  Since the centre of an Azumaya algebra is
$\mathcal O_S$, restriction gives a morphism of sheaves
\[
\mathrm{res}:\Der_F(\mathcal B)\longrightarrow\Der_F(\mathcal O_S).
\]
Its kernel is exactly $\Der_{\mathcal O_S}(\mathcal B)$, identified by the
first sequence with $\mathcal B/\mathcal O_S$. 
It remains to prove surjectivity as a sheaf statement.  Over a trivialising
\'{e}tale open $U$, every $F$-derivation $d$ of $\mathcal O_U$ has the entrywise
lift
\[
\widetilde d((a_{ij}))=(d(a_{ij})).
\]
This is an $F$-derivation of $M_n(\mathcal O_U)$ restricting to $d$ on
scalar matrices.  Thus $\mathrm{res}$ is locally surjective on the \'{e}tale
site, which is exactly surjectivity of the associated sheaf morphism.  On
overlaps two such local lifts differ by an $\mathcal O_S$-linear derivation
and hence by an inner derivation.  This identifies the kernel and completes the
second exact sequence.
\end{proof}

\subsection{Local splitting diagrams}
Let $T\to S$ be an étale cover with a trivialisation
\[
\theta:\mathcal B_T\xrightarrow{\sim}M_n(\mathcal O_T).
\]
The Azumaya condition is transported through
\[
\begin{tikzcd}[column sep=large,row sep=large]
\mathcal B_T\otimes\mathcal B_T^{\op}
\arrow[r,"\Phi_{\mathcal B_T}"] \arrow[d,"\theta\otimes\theta^{\op}"']
& \cEnd(\mathcal B_T) \arrow[d,"\cEnd(\theta)"]\\
M_n(\mathcal O_T)\otimes M_n(\mathcal O_T)^{\op}
\arrow[r,"\sim"']
& \cEnd(M_n(\mathcal O_T)).
\end{tikzcd}
\]
The lower map is the double-centraliser isomorphism.  The derivation sequence becomes
\[
\begin{tikzcd}[column sep=large]
0 \arrow[r]
& \mathfrak{pgl}_n(\mathcal O_T) \arrow[r,"\operatorname{ad}"]
& \Der_F(M_n(\mathcal O_T)) \arrow[r,"\mathrm{res}"]
& \Der_F(\mathcal O_T) \arrow[r] \arrow[l,bend left=35,"\delta\mapsto( a_{ij}\mapsto\delta(a_{ij}) )"]
& 0.
\end{tikzcd}
\]
On $T\times_ST$, two trivialisations differ by a
$\operatorname{PGL}_n$-cocycle, and the corresponding entrywise splittings
differ by an inner derivation.  Thus the exact sequence is intrinsic although
its local splittings are not.

Finite locally free generic families give classifying morphisms as follows.

\begin{Proposition}\label{prop:generic-classifying-map}
Let $R$ be a commutative $F$-algebra and let $\mathcal U=\bigoplus_g\mathcal U_g$ be a graded $R$-algebra satisfying $I$.  Suppose that on an open subset $U\subseteq\Spec R$, every $\mathcal U_g|_U$ is locally free of constant rank $d_g$.  Then $\mathcal U|_U$ defines an object of the groupoid $\PI^G_{I,\mathbf d}(U)$, equivalently a classifying morphism
\[
U\longrightarrow\PI^G_{I,\mathbf d}
\]
well-defined up to $2$-isomorphism.
If $\mathcal U|_U$ is Azumaya as an ungraded $\mathcal O_U$-algebra, the morphism factors through $\PI^{G,\Az}_{I,\mathbf d}$.
\end{Proposition}

\begin{proof}
This is Theorem~\ref{thm:PI-stack-classification} applied to the sheaf associated with $\mathcal U|_U$.  The second assertion follows from Theorem~\ref{thm:Azumaya-open}.
\end{proof}

\subsection{Brauer classes and local splittings}
\begin{Proposition}\label{prop:Azumaya-splitting-torsor}
Let $\mathcal B$ be an Azumaya $\mathcal O_S$-algebra of degree $n$.  The sheaf
$
\operatorname{Isom}_{\mathcal O_S\text{-alg}}
(M_n(\mathcal O_S),\mathcal B)
$
on the small \'{e}tale site is a $\operatorname{PGL}_n$-torsor.  Its image under the connecting map associated with
\[
1\longrightarrow\mathbb G_m\longrightarrow\operatorname{GL}_n
\longrightarrow\operatorname{PGL}_n\longrightarrow1
\]
is the Brauer class $[\mathcal B]\in H^2_{\acute et}(S,\mathbb G_m)$.  The class vanishes if and only if $\mathcal B\cong\cEnd(\mathcal P)$ for a rank-$n$ vector bundle $\mathcal P$ on $S$.
\end{Proposition}
\begin{proof}
Azumaya algebras are \'{e}tale-locally matrix algebras, and every automorphism of $M_n$ is inner, with automorphism group $\operatorname{PGL}_n$.  Hence the sheaf of trivialisations is a $\operatorname{PGL}_n$-torsor.  The obstruction to lifting this torsor through $\operatorname{GL}_n\to\operatorname{PGL}_n$ is the indicated boundary class.  A lift is the frame torsor of a rank-$n$ bundle $\mathcal P$, and descent identifies the associated endomorphism algebra with $\mathcal B$.
\end{proof}

\begin{Corollary}\label{cor:azumaya-tangent-decomposition}
If a chosen \'{e}tale trivialisation identifies $\mathcal B$ with $M_n(\mathcal O_S)$, then Theorem~\ref{thm:Azumaya-derivations} becomes
\[
0\longrightarrow\mathfrak{pgl}_n(\mathcal O_S)
\longrightarrow\Der_F(M_n(\mathcal O_S))
\longrightarrow\Der_F(\mathcal O_S)
\longrightarrow0,
\]
and the entrywise action of $\Der_F(\mathcal O_S)$ gives a splitting depending on the chosen trivialisation.
\end{Corollary}
\begin{proof}
The identification $M_n(\mathcal O_S)/\mathcal O_S\cong\mathfrak{pgl}_n(\mathcal O_S)$ gives the kernel.  Applying a derivation of $\mathcal O_S$ to matrix entries is a Lie-algebra section of the restriction map.
\end{proof}

\subsection{Orbit and form substacks}
Let $A$ be a finite-dimensional graded algebra of graded dimension $\mathbf d$, viewed as a point of $\AlgLaws^G_{\mathbf d}(I)$ for $I\subseteq\idG(A)$.  Let $\mathcal O_A=H_{\mathbf d}\cdot A$ be its change-of-basis orbit.

\begin{Proposition}\label{prop:orbit-classifying-stack}
The point $A:\Spec F\to\PI^G_{I,\mathbf d}$ induces a fully faithful morphism
\[
B\Aut_G(A)\longrightarrow\PI^G_{I,\mathbf d}.
\]
Its essential image over an $F$-scheme $S$ consists of the families which are fppf-locally isomorphic, as unital graded algebras, to $A\otimes_F\mathcal O_S$.  If the orbit fppf sheaf $\mathcal O_A=H_{\mathbf d}/\Aut_G(A)$ is represented, this morphism identifies its image with
\[
[\mathcal O_A/H_{\mathbf d}]\simeq B\Aut_G(A).
\]
\end{Proposition}
\begin{proof}
An $\Aut_G(A)$-torsor $P\to S$ produces by descent the graded algebra $P\times^{\Aut_G(A)}A$, which is fppf-locally the constant algebra $A$.  Conversely, the sheaf of graded algebra isomorphisms from $A\otimes_F\mathcal O_S$ to such a family is an $\Aut_G(A)$-torsor.  These constructions are quasi-inverse and identify automorphisms, proving full faithfulness and the essential-image statement.  The final assertion is the standard quotient description of a transitive orbit sheaf.
\end{proof}

\subsection{Central geometry and PI-degree}
\begin{Proposition}\label{prop:azumaya-PI-degree}
Let $\mathcal B$ be an Azumaya algebra of constant rank $n^2$ on an $F$-scheme $S$.  For every geometric point $\bar s\to S$ one has
\[
\mathcal B_{\bar s}\cong M_n(k(\bar s)),
\qquad
\operatorname{PIdeg}(\mathcal B_{\bar s})=n.
\]
If the residue characteristic is zero, the fibre satisfies the standard identity $s_{2n}$; for $n\geq2$ it does not satisfy $s_{2n-2}$.
\end{Proposition}
\begin{proof}
By the defining property of an Azumaya algebra, every geometric fibre is a central simple algebra of degree $n$ over an algebraically closed field and hence is $M_n$.  Its PI-degree is therefore $n$.  The final assertion is the Amitsur--Levitzki theorem together with the minimality of the standard identity for matrix algebras.
\end{proof}
\subsection{Brauer torsion and algebraic K-theory}
Here $\Br(S)$ denotes the Azumaya Brauer group of a scheme $S$.  For an Azumaya $\mathcal O_S$-algebra $\mathcal A$, let $\operatorname{Vect}(S,\mathcal A)$ denote the exact category of left $\mathcal A$-modules that are finite locally free as $\mathcal O_S$-modules.  Its conflations are the short exact sequences of $\mathcal A$-modules whose terms lie in this category.  Write
\[
\mathbf K(S;\mathcal A)=\mathbf K\bigl(\operatorname{Vect}(S,\mathcal A)\bigr),
\qquad K_i(S;\mathcal A)=\pi_i\mathbf K(S;\mathcal A),
\]
where $\mathbf K$ denotes connective Quillen K-theory; for $\mathcal A=\mathcal O_S$ we write $\mathbf K(S)$ and $K_i(S)$.

\begin{Theorem}\label{thm:brauer-torsion-on-PI-stack}
Let $N=\sum_{g\in G}d_g$ and set $\mathscr X=\PI^{G,\Az}_{I,\mathbf d}$.
\begin{enumerate}
\item If $N$ is not a square, then $\mathscr X(S)=\varnothing$ for every nonempty $F$-scheme $S$.
\item If $N=n^2$, every $\mathcal A\in\mathscr X(S)$ is an Azumaya algebra of degree $n$ and its Brauer class satisfies
$
n[\mathcal A]_{\mathrm{Br}}=0 \in \Br(S).
$
The maps
\[
\beta_{I,\mathbf d,S}:\pi_0\mathscr X(S)\longrightarrow\Br(S)[n],
\qquad [\mathcal A]_{\mathrm{iso}}\longmapsto[\mathcal A]_{\mathrm{Br}},
\]
are functorial in $S$.
\end{enumerate}
\end{Theorem}
\begin{proof}
Every object of $\PI^G_{I,\mathbf d}(S)$ is finite locally free of rank $N$ over $\mathcal O_S$.  If it is Azumaya, then \'{e}tale-locally it is isomorphic to $M_r(\mathcal O_S)$ and therefore has rank $r^2$.  This proves the first assertion and shows that $r=n$ when $N=n^2$.
Fix $\mathcal A\in\mathscr X(S)$ and choose an \'{e}tale cover $\{U_i\to S\}$ together with trivialisations $\phi_i:\mathcal A|_{U_i}\xrightarrow{\sim}M_n(\mathcal O_{U_i})$.  On $U_{ij}$ the automorphism $\phi_i\phi_j^{-1}$ defines an element $\alpha_{ij}\in\operatorname{PGL}_n(U_{ij})$.  After refining the cover, choose lifts $g_{ij}\in\operatorname{GL}_n(U_{ij})$.  On $U_{ijk}$ there are units $c_{ijk}$ such that
$
g_{ij}g_{jk}=c_{ijk}g_{ik}.
$
The \v{C}ech class of $(c_{ijk})$ is the image of the $\operatorname{PGL}_n$-torsor of trivialisations under the connecting map associated with
\[
1\longrightarrow\mathbb G_m\longrightarrow\operatorname{GL}_n
\longrightarrow\operatorname{PGL}_n\longrightarrow1,
\]
and hence represents $[\mathcal A]_{\mathrm{Br}}$.  Taking determinants gives
\[
c_{ijk}^{\,n}=\det(g_{ij})\det(g_{jk})\det(g_{ik})^{-1},
\]
so $(c_{ijk}^{\,n})$ is a coboundary.  Thus $n[\mathcal A]_{\mathrm{Br}}=0$; see also \cite{StacksBrauer}.  The torsor of trivialisations and the connecting class commute with pullback, which proves functoriality.  Isomorphic algebras determine isomorphic torsors, so the map factors through $\pi_0\mathscr X(S)$.
\end{proof}

\begin{Definition}\label{def:PI-Azumaya-Morita-groupoid}
Assume that $N=\sum_{g\in G}d_g=n^2$ and put $\mathscr X=\PI^{G,\Az}_{I,\mathbf d}$.  For an $F$-scheme $S$, let $\mathscr X^{\Mor}(S)$ be the Morita $(2,1)$-groupoid whose objects are the families $\mathcal A\in\mathscr X(S)$, whose $1$-morphisms $\mathcal A\to\mathcal A'$ are invertible $\mathcal A'$--$\mathcal A$-bimodules, and whose $2$-morphisms are bimodule isomorphisms.  The grading and the identities constrain the objects; Morita morphisms are taken between the underlying Azumaya algebras.  Let $\mathbf{BrAz}_n(S)$ be the analogous $(2,1)$-groupoid of all degree-$n$ Azumaya algebras on $S$.  There is a fully faithful inclusion
\[
\Phi_{I,\mathbf d,S}:\mathscr X^{\Mor}(S)\longrightarrow\mathbf{BrAz}_n(S).
\]
\end{Definition}

\begin{Theorem}\label{thm:PI-Azumaya-K-pseudofunctor}
For every $F$-scheme $S$, tensor product defines a pseudofunctor
\[
\mathbf{Vect}_S:\mathbf{BrAz}_n(S)\longrightarrow\mathbf{ExCat},
\qquad
\mathcal A\longmapsto\operatorname{Vect}(S,\mathcal A),
\qquad
{}_{\mathcal A'}\mathcal P_{\mathcal A}\longmapsto
\mathcal P\otimes_{\mathcal A}-,
\]
where $\mathbf{ExCat}$ denotes the $2$-category of exact categories, exact functors and natural isomorphisms.  Let $\mathbf{Sp}^{\simeq}$ denote the maximal $\infty$-groupoid of spectra.  Applying connective Quillen K-theory yields a homotopy-coherent functor
\[
\mathbf K_S:\mathbf{BrAz}_n(S)\longrightarrow\mathbf{Sp}^{\simeq},
\qquad
\mathcal A\longmapsto\mathbf K(S;\mathcal A).
\]
Its restriction to the PI-Azumaya moduli problem is the composite
\[
\mathbf K_{I,\mathbf d,S}:\mathscr X^{\Mor}(S)
\xrightarrow{\ \Phi_{I,\mathbf d,S}\ }
\mathbf{BrAz}_n(S)
\xrightarrow{\ \mathbf K_S\ }
\mathbf{Sp}^{\simeq}.
\]
In particular:
\begin{enumerate}
\item an invertible $\mathcal A'$--$\mathcal A$-bimodule induces an equivalence $\mathbf K(S;\mathcal A)\simeq\mathbf K(S;\mathcal A')$;
\item the map on connected components
\[
\pi_0\mathscr X^{\Mor}(S)\longrightarrow\Br(S),
\qquad[\mathcal A]_{\Mor}\longmapsto[\mathcal A]_{\mathrm{Br}},
\]
is injective and has image
$
\Br_{I,\mathbf d}(S):=
\bigl\{[\mathcal A]_{\mathrm{Br}}:\mathcal A\in\mathscr X(S)\bigr\}
\subseteq\Br(S)[n];
$
thus the equivalence class of $\mathbf K(S;\mathcal A)$ depends only on $[\mathcal A]_{\mathrm{Br}}\in\Br_{I,\mathbf d}(S)$;
\item for every morphism $f:T\to S$, pullback gives exact functors and maps of spectra
\[
f^*: \operatorname{Vect}(S,\mathcal A)\longrightarrow
\operatorname{Vect}(T,f^*\mathcal A),
\qquad
f^*: \mathbf K(S;\mathcal A)\longrightarrow
\mathbf K(T;f^*\mathcal A),
\]
and the canonical isomorphism
\[
f^*(\mathcal P\otimes_{\mathcal A}\mathcal M)
\cong f^*\mathcal P\otimes_{f^*\mathcal A}f^*\mathcal M
\]
makes the family $\mathbf K_{I,\mathbf d,S}$ compatible with base change, up to the canonical homotopies.
\end{enumerate}
\end{Theorem}
\begin{proof}
Let ${}_{\mathcal A'}\mathcal P_{\mathcal A}$ be invertible, with inverse ${}_{\mathcal A}\mathcal Q_{\mathcal A'}$.  The bimodule isomorphisms
$
\mathcal Q\otimes_{\mathcal A'}\mathcal P\cong\mathcal A,
$  $
\mathcal P\otimes_{\mathcal A}\mathcal Q\cong\mathcal A'
$
and the associativity constraint for tensor products give natural isomorphisms between the two composites and the corresponding identity functors.  On a common \'{e}tale splitting cover, write
\[
\mathcal A\cong\cEnd(\mathcal E),\qquad
\mathcal A'\cong\cEnd(\mathcal E'),\qquad
\mathcal P\cong\mathcal E'\otimes\mathcal E^{\vee}.
\]
If $e$ is a rank-one matrix idempotent on this cover, the usual Morita isomorphism identifies every $\mathcal M\in\operatorname{Vect}(S,\mathcal A)$ locally with $\mathcal E\otimes\mathcal W$, where $\mathcal W=e\mathcal M$ is a direct summand of the locally free $\mathcal O_S$-module $\mathcal M$ and is therefore locally free.  Hence
\[
\mathcal P\otimes_{\mathcal A}(\mathcal E\otimes\mathcal W)
\cong\mathcal E'\otimes\mathcal W.
\]
Thus $\mathcal P\otimes_{\mathcal A}-$ preserves finite local freeness.  It is exact because $\mathcal P$ is locally projective as a right $\mathcal A$-module.  The same argument applies to $\mathcal Q$.  For composable invertible bimodules $\mathcal P$ and $\mathcal R$, the canonical associator
\[
\mathcal R\otimes_{\mathcal A'}(\mathcal P\otimes_{\mathcal A}-)
\xrightarrow{\sim}
(\mathcal R\otimes_{\mathcal A'}\mathcal P)\otimes_{\mathcal A}-
\]
is the compositor; the usual unitors and the pentagon and triangle identities for tensor products prove the pseudofunctor axioms.  Quillen K-theory sends exact equivalences to equivalences of spectra and natural isomorphisms to homotopies, giving $\mathbf K_S$ and its restriction $\mathbf K_{I,\mathbf d,S}$ \cite{QuillenK}.

Two Azumaya algebras represent the same Brauer class if and only if they are Morita equivalent.  Hence two objects of $\mathscr X(S)$ lie in the same connected component of $\mathscr X^{\Mor}(S)$ precisely when their Brauer classes agree.  This proves injectivity on connected components.  Theorem~\ref{thm:brauer-torsion-on-PI-stack} places the image in $\Br(S)[n]$, and Morita invariance of $\mathbf K_S$ proves the second assertion.

Finally, a short exact sequence of finite locally free modules is locally split, so arbitrary pullback preserves the exact structures under consideration.  The displayed tensor--pullback isomorphism is functorial in the module and in the bimodule.  Its compatibility with associators and unitors is the standard coherence of tensor product; after applying K-theory it gives the asserted homotopies.  This proves base-change compatibility.
\end{proof}

\begin{Example}\label{ex:matrix-forms-K}
Assume that $G=\{e\}$, let $\mathbf d=(n^2)$ and take $I=\operatorname{Id}(M_n(F))$.  Then
\[
\mathscr X=\PI^{\Az}_{I,(n^2)}\simeq B\operatorname{PGL}_n.
\]
Indeed, for every $F$-scheme $S$, descent identifies an object of $B\operatorname{PGL}_n(S)$ with a $\operatorname{PGL}_n$-torsor $P\to S$, and the conjugation action of $\operatorname{PGL}_n$ on $M_n$ gives the associated algebra
\[
\mathcal A_P=P\times^{\operatorname{PGL}_n}M_n(\mathcal O_S).
\]
It is Azumaya of degree $n$ and satisfies every polynomial identity of $M_n(F)$, since both assertions may be checked fppf-locally, where $P$ is trivial and $\mathcal A_P\cong M_n(\mathcal O_S)$.  Conversely, an object of $\mathscr X(S)$ is an Azumaya algebra of degree $n$; its torsor of algebra trivialisations is a $\operatorname{PGL}_n$-torsor, and the two constructions are quasi-inverse.  Under this equivalence, Theorem~\ref{thm:PI-Azumaya-K-pseudofunctor} is the composite
\[
P\longmapsto\mathcal A_P\longmapsto\operatorname{Vect}(S,\mathcal A_P)\longmapsto\mathbf K(S;\mathcal A_P).
\]
The connecting morphism for $1\to\mathbb G_m\to\operatorname{GL}_n\to\operatorname{PGL}_n\to1$ sends $P$ to $[\mathcal A_P]_{\mathrm{Br}}$.  If $P$ is trivial, then $\mathcal A_P\cong M_n(\mathcal O_S)$, and the standard Morita bimodule $\mathcal O_S^n$ yields exact equivalences
\[
\operatorname{Vect}(S,M_n(\mathcal O_S))\simeq\operatorname{Vect}(S),
\qquad
\mathbf K(S;M_n(\mathcal O_S))\simeq\mathbf K(S).
\]
For a nontrivial torsor, the same construction gives the twisted spectrum $\mathbf K(S;\mathcal A_P)$; its equivalence class depends only on $[\mathcal A_P]_{\mathrm{Br}}$.
\end{Example}

\begin{Corollary}\label{cor:K-Azumaya-after-degree}
Let $S$ be connected and noetherian, and let $\mathcal A$ be an Azumaya algebra of degree $n$ on $S$.  Assume either that $S$ is regular or that $S$ has finite dimension and an ample invertible sheaf.  For every $i\geq0$, the exact functor $\mathcal E\mapsto\mathcal A\otimes_{\mathcal O_S}\mathcal E$ induces a homomorphism
\[
B_{\mathcal A}:K_i(S)\longrightarrow K_i(S;\mathcal A)
\]
whose kernel and cokernel are annihilated by a power of $n$.  Consequently,
\[
B_{\mathcal A}\otimes\mathbb Z[1/n]:
K_i(S)\otimes_{\mathbb Z}\mathbb Z[1/n]
\xrightarrow{\sim}
K_i(S;\mathcal A)\otimes_{\mathbb Z}\mathbb Z[1/n].
\]
\end{Corollary}
\begin{proof}
The $\mathcal O_S$-rank of $\mathcal A$ is $n^2$.  Hazrat--Hoobler prove, under the stated hypotheses, that the kernel and cokernel of $B_{\mathcal A}$ are annihilated by $(n^2)^m$ for some $m\geq0$ \cite{HazratHoobler}.  They are therefore annihilated by a power of $n$, and tensoring with $\mathbb Z[1/n]$ yields the asserted isomorphism.
\end{proof}

\section{Morita transport of PI-coefficients}\label{sec:morita-transport}

\begin{Definition}\label{def:Morita-2-groupoid}
A Morita $(2,1)$-groupoid $\mathscr M$ has unital $F$-algebras as selected objects, invertible bimodules ${}_{B'}P_B$ as $1$-morphisms $B\to B'$, and bimodule isomorphisms as $2$-morphisms.  Composition is
\[
({}_{B''}Q_{B'})\circ({}_{B'}P_B)
={}_{B''}(Q\otimes_{B'}P)_B,
\]
and the unit at $B$ is the regular bimodule ${}_BB_B$.
\end{Definition}

\begin{Proposition}\label{prop:Morita-module-transport}
Let ${}_{B'}P_B$ be invertible with inverse ${}_BQ_{B'}$.  Then
\[
P\otimes_B-:\Proj(B)^{\simeq}\xrightarrow{\sim}\Proj(B')^{\simeq}
\]
is an equivalence with quasi-inverse $Q\otimes_{B'}-$.  For every topological space $X$, sheaf extension of scalars
$
\underline P\otimes_{\underline B}-
$
preserves locally constant sheaves with finitely generated projective stalks, and its stalk at $x$ is canonically $P\otimes_B\mathcal E_x$.
\end{Proposition}
\begin{proof}
Choose bimodule isomorphisms
$
P\otimes_BQ\xrightarrow{\sim}B',
$ and $
Q\otimes_{B'}P\xrightarrow{\sim}B
$
compatible with the bimodule structures.  For a left $B$-module $M$,
associativity of tensor products gives natural isomorphisms
\[
Q\otimes_{B'}(P\otimes_BM)
\cong(Q\otimes_{B'}P)\otimes_BM
\cong B\otimes_BM\cong M.
\]
The analogous calculation for a left $B'$-module gives
$P\otimes_B(Q\otimes_{B'}N)\cong N$.  Hence the two tensor functors are
quasi-inverse equivalences.  Equivalences preserve projective objects and
finite generation; alternatively, an invertible bimodule is finitely
generated projective on the relevant sides and tensoring carries a finite
direct-summand presentation to one of the same kind.

Let $\mathcal E$ be a locally constant sheaf of left $\underline B$-modules.
On an open set $U$ on which
$\mathcal E|_U\cong\underline M|_U$, sheafified tensor product gives
\[
(\underline P\otimes_{\underline B}\mathcal E)|_U
\cong\underline{P\otimes_BM}|_U.
\]
Thus the result is locally constant and its stalks are finitely generated
projective.  For every $x\in X$, stalks commute with sheafification and with
tensor products of sheaves of modules, giving
\[
(\underline P\otimes_{\underline B}\mathcal E)_x
\cong P\otimes_B\mathcal E_x.
\]
The isomorphism is natural in $P$, $\mathcal E$, and $x$.
\end{proof}

Polynomial identities are not Morita invariant: $F$ and $M_n(F)$ are Morita equivalent but have different ordinary identities.  Hence $\mathscr M$ is a specified coefficient groupoid, not the Morita closure of a fixed PI-variety.

\section{Morita-equivariant monodromy over the coefficient base}\label{sec:morita-monodromy}

Let $X$ be connected, locally path-connected and semilocally simply connected, and let $\Pi_1(X)$ be its fundamental groupoid.  Let $\mathscr M$ be the Morita $(2,1)$-groupoid fixed in Section~\ref{sec:morita-transport}.  Its objects are unital $F$-algebras, its $1$-morphisms $B\to B'$ are invertible bimodules ${}_{B'}P_B$, and its $2$-morphisms are bimodule isomorphisms.

For an object $B$ of $\mathscr M$, let $\Proj(B)^{\simeq}$ be the core groupoid of finitely generated projective left $B$-modules.  Let $\LCS(X,B)$ be the groupoid of locally constant sheaves of left modules over the constant sheaf $\underline B$ whose stalks are finitely generated projective.

\begin{Proposition}\label{prop:Morita-pseudofunctors}
The assignments
\[
\mathcal L_X(B)=\LCS(X,B),\qquad
\mathcal M_X(B)=\Fun(\Pi_1(X),\Proj(B)^{\simeq})
\]
and, for ${}_{B'}P_B$,
\[
\mathcal L_X(P)(\mathcal E)=\underline P\otimes_{\underline B}\mathcal E,
\qquad
\mathcal M_X(P)(F)=(P\otimes_B-)\circ F,
\]
define pseudofunctors $\mathcal L_X,\mathcal M_X:\mathscr M\to\Gpd$.  Their compositor at composable bimodules $P,Q$ is induced by
\[
\underline Q\otimes_{\underline{B'}}
(\underline P\otimes_{\underline B}-)
\xrightarrow{\sim}
\underline{Q\otimes_{B'}P}\otimes_{\underline B}-
\]
on the sheaf side and by the analogous associativity isomorphism on projective modules; the unitor is induced by $B\otimes_B-\cong\id$.
\end{Proposition}
\begin{proof}
Proposition~\ref{prop:Morita-module-transport} shows that every displayed
extension-of-scalars functor preserves the specified objects and is an
equivalence of groupoids.  If
$\theta:P\xrightarrow{\sim}P'$ is a bimodule isomorphism, the maps
\[
\theta\otimes\id:
P\otimes_BM\xrightarrow{\sim}P'\otimes_BM
\]
form a natural isomorphism, and the same formula after sheafification gives
the action on $2$-morphisms.
For composable arrows $B\xrightarrow{P}B'\xrightarrow{Q}B''$, define the
compositor on a sheaf $\mathcal E$ by
\[
q\otimes(p\otimes e)\longmapsto(q\otimes p)\otimes e.
\]
This map is balanced over $B'$ and $B$, hence descends to the displayed
isomorphism of sheaf tensor products.  The same formula defines the
compositor on projective modules.  For three bimodules $P,Q,R$, both paths
around the pentagon send
$r\otimes(q\otimes(p\otimes m))$ to
$((r\otimes q)\otimes p)\otimes m$ under the canonical associativity
identifications.  Hence the pentagon commutes.  The unitors are the maps
$b\otimes m\mapsto bm$, and the two triangle diagrams commute on pure
tensors.  Naturality then gives all coherence axioms.  Therefore the
assignments define pseudofunctors $\mathscr M\to\Gpd$.
\end{proof}

\subsection{Coherence diagrams for Morita composition}
For composable Morita arrows $B\xrightarrow{P}B'\xrightarrow{Q}B''$, the compositor of $\mathcal L_X$ is
\[
a_{Q,P,\mathcal E}:
\underline Q\otimes_{\underline{B'}}(\underline P\otimes_{\underline B}\mathcal E)
\xrightarrow{\sim}
\underline{Q\otimes_{B'}P}\otimes_{\underline B}\mathcal E.
\]
For three composable bimodules, coherence is the commutativity of
\[
\begin{tikzcd}[column sep=huge,row sep=large]
\underline R\otimes(\underline Q\otimes(\underline P\otimes\mathcal E))
\arrow[r,"a_{R,Q}"] \arrow[d,"\id\otimes a_{Q,P}"']
& \underline{R\otimes Q}\otimes(\underline P\otimes\mathcal E)
\arrow[d,"a_{R\otimes Q,P}"]\\
\underline R\otimes(\underline{Q\otimes P}\otimes\mathcal E)
\arrow[r,"a_{R,Q\otimes P}"']
& \underline{R\otimes Q\otimes P}\otimes\mathcal E.
\end{tikzcd}
\]
The pseudonaturality cells $\Theta_P$ satisfy
\[
\begin{tikzcd}[column sep=4.2em,row sep=large]
\operatorname{Mon}_{B''}\!\bigl(\underline Q\otimes(\underline P\otimes\mathcal E)\bigr)
\arrow[r,"\Theta_Q"] \arrow[d,"\operatorname{Mon}(a_{Q,P})"']
& (Q\otimes_{B'}-)\operatorname{Mon}_{B'}(\underline P\otimes\mathcal E)
\arrow[d,"Q\otimes\Theta_P"]\\
\operatorname{Mon}_{B''}\!\bigl(\underline{Q\otimes P}\otimes\mathcal E\bigr)
\arrow[r,"\Theta_{Q\otimes P}"']
& ((Q\otimes_{B'}P)\otimes_B-)\operatorname{Mon}_B(\mathcal E).
\end{tikzcd}
\]
Here the right vertical arrow includes the standard associator
$(Q\otimes_{B'}-)(P\otimes_B-)\cong((Q\otimes_{B'}P)\otimes_B-)$.
The diagram shows that the fibrewise equivalences are compatible over the Morita base.

\subsection{The two Grothendieck totals}
The bicategorical Grothendieck total $\int\mathcal L_X$ has objects $(B,\mathcal E)$ with $B\in\mathscr M$ and $\mathcal E\in\LCS(X,B)$.  A $1$-morphism
\[
(B,\mathcal E)\longrightarrow(B',\mathcal E')
\]
is a pair $(P,\alpha)$ consisting of a Morita bimodule ${}_{B'}P_B$ and an isomorphism
\begin{equation}\label{eq:intLX-1mor-expanded}
\alpha:\underline P\otimes_{\underline B}\mathcal E\xrightarrow{\sim}\mathcal E'.
\end{equation}
A $2$-morphism $(P,\alpha)\Rightarrow(P',\alpha')$ is a bimodule isomorphism $\theta:P\to P'$ satisfying
\begin{equation}\label{eq:intLX-2mor-expanded}
\alpha'\circ(\underline\theta\otimes\id_{\mathcal E})=\alpha.
\end{equation}
If $(Q,\beta):(B',\mathcal E')\to(B'',\mathcal E'')$, the composite is
\[
\bigl(Q\otimes_{B'}P,\ \beta\circ(\underline Q\otimes\alpha)\circ\phi_{Q,P}(\mathcal E)^{-1}\bigr),
\]
where
$\phi_{Q,P}:\mathcal L_X(Q)\mathcal L_X(P)\xrightarrow{\sim}\mathcal L_X(Q\otimes_{B'}P)$
is the compositor defined above.

The total $\int\mathcal M_X$ is defined analogously.  Its objects are pairs $(B,F)$ with $F:\Pi_1(X)\to\Proj(B)^{\simeq}$.  A $1$-morphism $(B,F)\to(B',F')$ is $(P,\alpha)$ with
\begin{equation}\label{eq:intMX-1mor-expanded}
\alpha:(P\otimes_B-)\circ F\xRightarrow{\sim}F'.
\end{equation}
The $2$-morphisms are bimodule isomorphisms compatible with these natural transformations, and composition is governed by
\[
Q\otimes_{B'}(P\otimes_B-)\cong(Q\otimes_{B'}P)\otimes_B-.
\]
Both totals project to the coefficient base:
\[
\pi_{\mathcal L}:\int\mathcal L_X\to\mathscr M,
\qquad
\pi_{\mathcal M}:\int\mathcal M_X\to\mathscr M.
\]
The projections fit into the diagram
\[
\begin{tikzcd}[column sep=huge,row sep=large]
\int\mathcal L_X \arrow[dr,"\pi_{\mathcal L}"'] \arrow[rr,dashed,"\mathbf{Mon}_{X;\mathscr M}"]
&&\int\mathcal M_X\arrow[dl,"\pi_{\mathcal M}"]\\
&\mathscr M&
\end{tikzcd}
\]
A morphism in either total records both the coefficient transport and an identification in the target fibre.  In particular, $(P,\alpha)$ factors through the canonical coCartesian arrow $(P,\id)$ followed by the vertical arrow induced by $\alpha$.

\subsection{The fixed-coefficient equivalence}
\begin{Theorem}\label{thm:fixed-coefficient-monodromy}
Let $X$ be locally path-connected and semilocally simply connected.  For every unital algebra $B$, parallel transport defines an equivalence of groupoids
\[
\operatorname{Mon}_B:\LCS(X,B)
\xrightarrow{\sim}
\Fun(\Pi_1(X),\Proj(B)^{\simeq}).
\]
A quasi-inverse sends $F$ to the sheaf of path-compatible families; on a path-connected open set $U$ whose loops are null-homotopic in $X$, evaluation at any base point identifies its sections with the corresponding value of $F$.
\end{Theorem}
\begin{proof}
Let $\mathcal E\in\LCS(X,B)$.  Choose a chain of connected trivialising
open sets along a path $\gamma:[0,1]\to X$.  The constant identifications on
successive overlaps give an isomorphism
$T_{\mathcal E}(\gamma):\mathcal E_{\gamma(0)}\to
\mathcal E_{\gamma(1)}$.  Refining the subdivision does not change this
map.  A homotopy of paths can be subdivided into rectangles contained in
trivialising opens, and transport around the boundary of each rectangle is
the identity; hence the map depends only on the homotopy class relative to
the endpoints.  Concatenation of paths becomes composition.  This defines a
functor
\[
\operatorname{Mon}_B(\mathcal E):\Pi_1(X)\to\Proj(B)^{\simeq}.
\]
A morphism of local systems commutes with all continuation maps, so this
construction is functorial.

Conversely, let $F:\Pi_1(X)\to\Proj(B)^{\simeq}$.  Define a presheaf
$\mathcal E_F$ by letting $\mathcal E_F(U)$ be the set of families
$s=(s_x)_{x\in U}$ with $s_x\in F(x)$ such that
\[
F([\gamma])(s_x)=s_y
\]
for every path $\gamma$ in $U$ from $x$ to $y$.  The condition is stable
under restriction, and compatible families glue pointwise, so
$\mathcal E_F$ is a sheaf of left $\underline B$-modules.
Because $X$ is locally path-connected and semilocally simply connected,
every point has a basis of path-connected open neighbourhoods $U$ for which
every loop in $U$ is null-homotopic in $X$.  On such a $U$, choose
$x_0\in U$.  Evaluation at $x_0$ gives an isomorphism
\[
\mathcal E_F(U)\xrightarrow{\sim}F(x_0):
\]
a vector $v\in F(x_0)$ determines the unique section
$s_x=F([\gamma_{x_0,x}])(v)$, because the value is independent of the path
inside $U$.  Thus $\mathcal E_F$ is locally constant with stalk
$(\mathcal E_F)_x\cong F(x)$, and its stalks are finitely generated
projective.
Under these identifications, the continuation map of $\mathcal E_F$ along a
path $\gamma$ is exactly $F([\gamma])$.  Hence
$\operatorname{Mon}_B(\mathcal E_F)\cong F$.  Conversely, a local section
of a local system is exactly a family of germs invariant under its local
parallel transport, which gives a canonical isomorphism
$\mathcal E_{\operatorname{Mon}_B(\mathcal E)}\cong\mathcal E$.
The constructions on morphisms are also inverse: a natural transformation
is evaluated pointwise on stalks and glues on the monodromy-trivial opens.
Therefore $\operatorname{Mon}_B$ is an equivalence of groupoids.
\end{proof}

\subsection{CoCartesian factorisation}
A morphism $(P,\alpha):(B,\mathcal E)\to(B',\mathcal E')$ in $\int\mathcal L_X$ factors canonically as
\[
(B,\mathcal E)\xrightarrow{(P,\id)}
(B',\underline P\otimes_{\underline B}\mathcal E)
\xrightarrow{(B',\alpha)}(B',\mathcal E').
\]
The first arrow is coCartesian over $P$ and the second is vertical.  The analogous factorisation on the Betti side is
\[
(B,F)\xrightarrow{(P,\id)}
(B',(P\otimes_B-)\circ F)
\xrightarrow{(B',\alpha)}(B',F').
\]
The comparison cell of Lemma~\ref{lem:expanded-stalk-morita} carries the first factorisation to the second.  Thus the biequivalence lies over the Morita base, rather than only after the projections are forgotten.

\subsection{Description of hom-groupoids}
Fix $(B,\mathcal E)$ and $(B',\mathcal E')$.  The hom-groupoid in $\int\mathcal L_X$ is the Grothendieck groupoid over the hom-groupoid $\mathscr M(B,B')$: its objects over a Morita arrow $P:B\to B'$ are the isomorphisms
\[
\underline P\otimes_{\underline B}\mathcal E\xrightarrow{\sim}\mathcal E',
\]
with morphisms induced by bimodule isomorphisms.  The corresponding hom-groupoid on the Betti side uses
\[
(P\otimes_B-)\circ\operatorname{Mon}_B(\mathcal E)
\xRightarrow{\sim}
\operatorname{Mon}_{B'}(\mathcal E').
\]
For each fixed $P$, the pseudonaturality cell identifies these two groupoids through the full faithfulness of $\operatorname{Mon}_{B'}$.  This fibrewise description proves local equivalence of the total bicategories.

\subsection{Compatibility of stalks and transport}
\begin{Lemma}\label{lem:expanded-stalk-morita}
For $\mathcal E\in\LCS(X,B)$ and $x\in X$, there is a canonical natural isomorphism
\[
\Theta_x:(\underline P\otimes_{\underline B}\mathcal E)_x\xrightarrow{\sim}P\otimes_B\mathcal E_x.
\]
For every path class $[\gamma]:x\to y$, the square
\[
\begin{CD}
(\underline P\otimes_{\underline B}\mathcal E)_x
@>{T_{\underline P\otimes\mathcal E}([\gamma])}>>
(\underline P\otimes_{\underline B}\mathcal E)_y\\
@V{\Theta_x}VV @VV{\Theta_y}V\\
P\otimes_B\mathcal E_x
@>{\id_P\otimes T_{\mathcal E}([\gamma])}>>
P\otimes_B\mathcal E_y
\end{CD}
\]
commutes.
\end{Lemma}
\begin{proof}
The stalk formula follows because the sheaf tensor product is obtained from the presheaf tensor product by sheafification, stalks are unchanged by sheafification, and $P\otimes_B-$ commutes with filtered colimits.  For transport, subdivide a representative path into pieces contained in open sets on which $\mathcal E$ is constant.  On each piece both local systems are constant and the comparison is the identity on $P\otimes_BM$.  Pasting gives the square for the full path.
\end{proof}

The lemma gives a natural isomorphism
\begin{equation}\label{eq:expanded-pseudonaturality-cell}
\operatorname{Mon}_{B'}(\underline P\otimes_{\underline B}\mathcal E)
\xRightarrow{\sim}
(P\otimes_B-)\circ\operatorname{Mon}_B(\mathcal E).
\end{equation}
It is natural in $\mathcal E$, compatible with bimodule isomorphisms and coherent under tensor products of Morita bimodules.

\begin{Theorem}\label{thm:expanded-Morita-monodromy}
The fixed-coefficient monodromy equivalences
$
\operatorname{Mon}_B:\LCS(X,B)\xrightarrow{\sim}
\Fun(\Pi_1(X),\Proj(B)^{\simeq})
$
assemble, with the cells \eqref{eq:expanded-pseudonaturality-cell}, into a pseudonatural equivalence
$
\operatorname{Mon}:\mathcal L_X\Longrightarrow\mathcal M_X.
$
They induce a biequivalence over $\mathscr M$
\[
\mathbf{Mon}_{X;\mathscr M}:\int\mathcal L_X\xrightarrow{\sim}\int\mathcal M_X,
\qquad
\pi_{\mathcal M}\circ\mathbf{Mon}_{X;\mathscr M}=\pi_{\mathcal L}.
\]
Its restriction to the fibre over $B$ is $\operatorname{Mon}_B$.
\end{Theorem}
\begin{proof}
For each coefficient algebra $B$, Theorem~\ref{thm:fixed-coefficient-monodromy}
provides an equivalence $\operatorname{Mon}_B$.  Let
${}_{B'}P_B$ be a Morita arrow and $\mathcal E\in\LCS(X,B)$.  At a point
$x$, the canonical stalk isomorphism is
\[
\Theta_{P,\mathcal E,x}:
(\underline P\otimes_{\underline B}\mathcal E)_x
\xrightarrow{\sim}P\otimes_B\mathcal E_x.
\]
For a path $[\gamma]:x\to y$, local triviality shows that transport on the
left sends the germ of $p\otimes e$ to the germ of
$p\otimes T_{\mathcal E}([\gamma])e$.  Therefore
\[
\Theta_{P,\mathcal E,y}\circ
T_{\underline P\otimes\mathcal E}([\gamma])
=(\id_P\otimes T_{\mathcal E}([\gamma]))
\circ\Theta_{P,\mathcal E,x}.
\]
The $\Theta_{P,\mathcal E,x}$ thus form the invertible pseudonaturality cell
\[
\Theta_P:\operatorname{Mon}_{B'}\circ\mathcal L_X(P)
\Longrightarrow\mathcal M_X(P)\circ\operatorname{Mon}_B.
\]
They are natural in $\mathcal E$ and in bimodule isomorphisms.

For composable $P$ and $Q$, the two composites required by pseudonaturality
agree on every stalk and on every pure tensor $q\otimes p\otimes e$; both
send it to $(q\otimes p)\otimes e$.  The identity coherence is equally the
map $b\otimes e\mapsto be$.  Hence the family
$(\operatorname{Mon}_B,\Theta_P)$ is a pseudonatural equivalence.
The induced homomorphism on Grothendieck totals is as follows.  On objects
put
\[
(B,\mathcal E)\longmapsto
(B,\operatorname{Mon}_B(\mathcal E)).
\]
A $1$-morphism $(P,\alpha)$ is sent to the composite natural isomorphism
\[
(P\otimes_B-)\operatorname{Mon}_B(\mathcal E)
\xRightarrow{\ \Theta_P^{-1}\ }
\operatorname{Mon}_{B'}(\underline P\otimes_{\underline B}\mathcal E)
\xRightarrow{\ \operatorname{Mon}_{B'}(\alpha)\ }
\operatorname{Mon}_{B'}(\mathcal E').
\]
On $2$-morphisms it applies the natural transformation induced by the
bimodule isomorphism.  The coherence established above proves compatibility with
composition and units, and the construction lies over the identity of
$\mathscr M$.

Essential surjectivity follows fibrewise: for $(B,F)$ choose
$\mathcal E\in\LCS(X,B)$ and an isomorphism
$\operatorname{Mon}_B(\mathcal E)\cong F$.  To prove local equivalence, fix
$(B,\mathcal E)$ and $(B',\mathcal E')$.  For each Morita arrow
$P:B\to B'$, the $P$-component of the source hom-groupoid is the groupoid of
isomorphisms
$\underline P\otimes\mathcal E\xrightarrow{\sim}\mathcal E'$, while the
corresponding target component is the groupoid of natural isomorphisms
$(P\otimes_B-)\operatorname{Mon}_B(\mathcal E)\to
\operatorname{Mon}_{B'}(\mathcal E')$.  The cell $\Theta_P$ identifies the
second with the groupoid of arrows between the monodromies of
$\underline P\otimes\mathcal E$ and $\mathcal E'$.  Full faithfulness of
$\operatorname{Mon}_{B'}$ identifies this with the first groupoid, including
its morphisms induced by bimodule isomorphisms.  Hence every hom-groupoid map
is an equivalence.  The induced homomorphism of totals is therefore a
biequivalence over $\mathscr M$, and its fibre over $B$ is
$\operatorname{Mon}_B$.
\end{proof}

\begin{Corollary}\label{cor:expanded-Morita-family}
The projections $\pi_{\mathcal L}$ and $\pi_{\mathcal M}$ are bicategorical opfibrations in groupoids.  Their fibres over $B$ are respectively $\LCS(X,B)$ and $\Fun(\Pi_1(X),\Proj(B)^{\simeq})$.  CoCartesian transport along ${}_{B'}P_B$ is extension of scalars on the sheaf side and postcomposition with $P\otimes_B-$ on the Betti side.  The biequivalence of Theorem~\ref{thm:expanded-Morita-monodromy} identifies these opfibration structures.
\end{Corollary}
\begin{proof}
This is the Grothendieck construction for groupoid-valued pseudofunctors; see \cite{Benabou,StreetFib,VistoliNotes}.  The canonical arrows $(P,\id)$ are the chosen coCartesian lifts.  A biequivalence over the base preserves and reflects the coCartesian universal property.
\end{proof}

\subsection{Exact monodromy and K-theory}
Let $\LCS^{\mathrm{proj}}_{\mathrm{ex}}(X,B)$ be the additive category of locally constant sheaves of left $\underline B$-modules with finitely generated projective stalks, with all sheaf morphisms.  Declare a kernel--cokernel pair to be a conflation when it is exact on every stalk.  Equip $\Proj(B)$ with its split exact structure and $\Fun(\Pi_1(X),\Proj(B))$ with the objectwise exact structure.

\begin{Theorem}\label{thm:K-monodromy}
The monodromy construction extends to an exact equivalence
\[
\operatorname{Mon}^{\mathrm{ex}}_B:\LCS^{\mathrm{proj}}_{\mathrm{ex}}(X,B)
\xrightarrow{\sim}\Fun(\Pi_1(X),\Proj(B)).
\]
In particular, the stalkwise conflations define an exact structure on $\LCS^{\mathrm{proj}}_{\mathrm{ex}}(X,B)$, and monodromy induces an equivalence of connective K-theory spectra
\[
\mathbf K\bigl(\LCS^{\mathrm{proj}}_{\mathrm{ex}}(X,B)\bigr)
\simeq
\mathbf K\bigl(\Fun(\Pi_1(X),\Proj(B))\bigr).
\]
For every Morita arrow ${}_{B'}P_B$, the natural comparison between monodromy and extension of scalars induces a homotopy-commutative square
\[
\begin{tikzcd}[column sep=large,row sep=large]
\mathbf K\bigl(\LCS^{\mathrm{proj}}_{\mathrm{ex}}(X,B)\bigr)
\arrow[r,"\mathbf K(\operatorname{Mon}^{\mathrm{ex}}_B)"]
\arrow[d,"\mathbf K(\underline P\otimes_{\underline B}-)"']
& \mathbf K\bigl(\Fun(\Pi_1(X),\Proj(B))\bigr)
\arrow[d,"\mathbf K(P\otimes_B-)" ]\\
\mathbf K\bigl(\LCS^{\mathrm{proj}}_{\mathrm{ex}}(X,B')\bigr)
\arrow[r,"\mathbf K(\operatorname{Mon}^{\mathrm{ex}}_{B'})"']
& \mathbf K\bigl(\Fun(\Pi_1(X),\Proj(B'))\bigr).
\end{tikzcd}
\]
Thus the induced isomorphisms on $K_i$, $i\geq0$, are compatible with Morita transport.
\end{Theorem}
\begin{proof}
The proof of Theorem~\ref{thm:fixed-coefficient-monodromy} applies without restricting morphisms to isomorphisms.  A sheaf morphism is sent to its maps on stalks, and a natural transformation between fundamental-groupoid functors glues uniquely on monodromy-trivial open sets.  These operations are mutually inverse and give an equivalence of additive categories.
Under this equivalence, a sequence of local systems is exact on every stalk if and only if the corresponding sequence of functors is exact at every object of $\Pi_1(X)$.  The latter sequences form the standard objectwise exact structure on the functor category.  Transporting that exact structure across the equivalence proves the first assertion and shows that $\operatorname{Mon}^{\mathrm{ex}}_B$ is exact.
The stalk isomorphisms of Lemma~\ref{lem:expanded-stalk-morita} are natural for arbitrary morphisms and therefore define a natural isomorphism of exact functors
\[
\operatorname{Mon}^{\mathrm{ex}}_{B'}\circ
(\underline P\otimes_{\underline B}-)
\xRightarrow{\sim}
(P\otimes_B-)\circ\operatorname{Mon}^{\mathrm{ex}}_B.
\]
Quillen K-theory sends the exact equivalence to an equivalence of spectra and this natural isomorphism to the homotopy displayed above \cite{QuillenK}.
\end{proof}

\begin{Proposition}\label{prop:base-change-of-monodromy}
Let $\mathscr N\to\mathscr M$ be a homomorphism of Morita $(2,1)$-groupoids.  Pullback of the pseudofunctors $\mathcal L_X$ and $\mathcal M_X$ to $\mathscr N$ carries the pseudonatural monodromy equivalence to the corresponding equivalence over $\mathscr N$.  Hence the biequivalence of totals is stable under restriction to any Morita subgroupoid, in particular to a subgroupoid generated by selected Azumaya fibres or by representatives of a fixed Brauer class.
\end{Proposition}
\begin{proof}
Pseudonatural transformations and their coherence cells pull back by precomposition.  The Grothendieck construction is functorial for such base change, and the fibrewise equivalences remain equivalences.
\end{proof}

\subsection{Quantum-torus specialisation}
\begin{Proposition}\label{prop:quantum-torus-Azumaya-Morita}
Let $q\in\mathbb C^*$ be a primitive $n$-th root of unity and
\[
B_q=\mathbb C\langle U^{\pm1},V^{\pm1}\rangle/(UV-qVU),
\qquad a=U^n,\quad b=V^n,\quad R=\mathbb C[a^{\pm1},b^{\pm1}].
\]
Then $Z(B_q)=R$, and $B_q$ is an Azumaya $R$-algebra of rank $n^2$.  The extension
$
R'=R[u^{\pm1},v^{\pm1}]/(u^n-a,v^n-b)
$
is finite \'{e}tale and faithfully flat, and there is an $R'$-algebra isomorphism
\[
B_q\otimes_RR'\xrightarrow{\sim}M_n(R').
\]
Consequently, for every connected, locally path-connected and semilocally simply connected space $X$, the standard Morita bimodule induces an equivalence between the Betti groupoids with coefficients in $B_q\otimes_RR'$ and in $R'$.
\end{Proposition}
\begin{proof}
The monomials $U^iV^j$, $0\leq i,j<n$, form an $R$-basis of $B_q$.  If an $R$-linear combination of these monomials commutes with both $U$ and $V$, primitivity of $q$ forces $i=j=0$ for every nonzero coefficient; hence $Z(B_q)=R$.

The elements $u$ and $v$ are units, and the derivatives of $T^n-a$ and $T^n-b$ are invertible at their roots.  Thus $R'/R$ is finite \'{e}tale and faithfully flat.  Let
\[
D=\operatorname{diag}(1,q,\ldots,q^{n-1}),
\qquad Se_j=e_{j+1}\quad(j\bmod n).
\]
Then $DS=qSD$ and $D^n=S^n=I_n$.  Therefore
\[
U\longmapsto uD,\qquad V\longmapsto vS
\]
defines an $R'$-algebra homomorphism $B_q\otimes_RR'\to M_n(R')$.  The matrices $D$ and $S$ generate $M_n(\mathbb C)$, so the homomorphism is surjective.  Both sides are finite locally free of rank $n^2$ over $R'$, hence it is an isomorphism.  Faithfully flat descent now shows that $B_q$ is Azumaya over $R$.  The final assertion follows from the standard Morita equivalence between $M_n(R')$ and $R'$ together with Theorem~\ref{thm:expanded-Morita-monodromy}.
\end{proof}

\begin{Corollary}\label{cor:quantum-torus-K}
For every $i\geq0$,
\[
K_i(\Spec R)\otimes_{\mathbb Z}\mathbb Z[1/n]
\xrightarrow{\sim}
K_i(\Spec R;B_q)\otimes_{\mathbb Z}\mathbb Z[1/n].
\]
After the splitting extension $R\to R'$, Morita equivalence gives an integral isomorphism $K_i(\Spec R';B_q\otimes_RR')\cong K_i(\Spec R')$.
\end{Corollary}
\begin{proof}
The Laurent polynomial ring $R$ is connected, noetherian and regular, and $B_q$ has degree $n$; the first assertion is Corollary~\ref{cor:K-Azumaya-after-degree}.  The second follows from $B_q\otimes_RR'\cong M_n(R')$ and Theorem~\ref{thm:PI-Azumaya-K-pseudofunctor}.
\end{proof}

For an elliptic curve $E$ and $X=E(\mathbb C)$, an object of
$
\Fun(\Pi_1(X),\Proj(R')^{\simeq})
$
is a finite projective $R'$-module together with two commuting automorphisms.  Morita transport sends
\[
(M,g_1,g_2)\longmapsto(R'^n\otimes_{R'}M,\id\otimes g_1,\id\otimes g_2).
\]
At a central specialisation, a projective $M_n(\mathbb C)$-module has the form $\mathbb C^n\otimes W$ and has automorphism group $\GL(W)$; the Betti rank is therefore $\dim W$, not automatically $n$.

\section{Examples and applications}\label{sec:examples}

The examples below illustrate closed PI-substacks, Azumaya and form loci, relatively free sheaves and filtered objects outside fixed finite rank.

\subsection{Classical and supergeometric structures}
For the trivial grading and $I=([x,y])$, $\PI^1_{I,d}$ parametrises commutative algebra structures on rank-$d$ bundles.  Structure sheaves of schemes or manifolds usually have infinite $F$-rank and belong to the sheaf-theoretic theory, whereas finite locally free commutative algebra schemes define points of the stack.  For $G=\mathbb Z_2$, a locally free supercommutative algebra of superrank $(p,q)$ gives a map to the corresponding closed PI-stack.  Split supermanifolds are locally modelled by exterior algebras and nonsplit ones by descent; for $\mathbb Z_2^n$, a sign bicharacter gives the colour-commutative equations.

\subsection{Matrix, Azumaya and Brauer geometry}
The identity ideal of $M_n(F)$ defines a closed PI-substack.  Its Azumaya open classifies families which are locally matrix algebras of degree $n$.  The centre supplies the ordinary local geometry, the algebra supplies the Brauer twist, and the orbit substack of the split matrix algebra is $B\operatorname{PGL}_n$.  Nontrivial forms are $\operatorname{PGL}_n$-torsors and may represent nonzero Brauer classes.  Universal differentials and the derivation sequence separate inner matrix directions from directions of the centre.

\subsection{Upper triangular and radical-containing algebras}
The algebra $UT_n(F)$ is finite-dimensional and PI\@.  Its moduli point is not Azumaya and has a nonzero radical.  The surrounding PI-stack records deformations that preserve the selected identity theory, while its orbit and normal tangent space measure rigidity and degeneration.

\subsection{Clifford and quaternion families}
Let $\mathcal E$ be locally free of rank $r$ and let $q:\operatorname{Sym}^2\mathcal E\to\mathcal O_X$ be a quadratic form.  The Clifford sheaf
\[
\operatorname{Cl}(\mathcal E,q)=T(\mathcal E)/(v\otimes v-q(v))
\]
is locally free of rank $2^r$, and its even part has rank $2^{r-1}$.  On the appropriate nondegenerate locus the even Clifford algebra is Azumaya over its centre.  In particular, for a nondegenerate rank-$3$ quadratic bundle, the even Clifford algebra has rank $4$ and is a quaternion Azumaya algebra.  The full Clifford algebra of a rank-$4$ bundle has rank $16$ and should not be identified with a quaternion algebra.  These sheaves yield classifying maps to finite-rank PI-stacks on loci of constant rank and provide basic examples of the Azumaya geometry \cite{Chan-Ingalls,HS,Lieb,Reede}.

\subsection{Grassmann and Fedosov models}
The infinite Grassmann algebra is a fundamental relatively free model, but its infinite rank places it in the sheaf-theoretic rather than the fixed-rank theory.  Likewise, let
\[
\Omega(\mathbb C^n)=\bigwedge(\mathbb C^n)\otimes\mathbb C[x_1,\ldots,x_n]
\]
and let $\Omega^{\mathrm{even}}$ be the even part.  Feigin and Shoikhet define
\[
\alpha*\beta=\alpha\wedge\beta+\frac12d\alpha\wedge d\beta,
\qquad [\alpha,\beta]_* = d\alpha\wedge d\beta.
\]
For $A_n=F\langle x_1,\ldots,x_n\rangle$, their map identifies the relevant quotient by the third lower-central-series term with this Fedosov algebra of even forms \cite{FeSh,AEt}.  The model is filtered and infinite-dimensional.  Assigning total degree $1$ to both $x_i$ and $dx_i$ makes the product homogeneous, so quotienting by terms of degree $>N$ gives finite-dimensional algebra quotients; exterior-degree truncation alone does not.

\subsection{Quantum tori and quantum groups at roots of unity}
Quantum tori at roots of unity are finite over their Laurent centres and yield explicit Azumaya families, linking PI-degree, central specialisation and Morita transport.  For quantum groups, three mechanisms must be distinguished:
\begin{enumerate}
\item the small quantum group is finite-dimensional and hence PI;
\item certain unrestricted or De Concini--Kac forms are finite over large central subalgebras;
\item quantum coordinate algebras and quantum affine spaces often contain central powers of the generators and are finite over the resulting centre.
\end{enumerate}
The latter two mechanisms produce central families whose generic fibres may be central simple and whose Azumaya loci may be studied by localisation.  They have different centres and representation theories, and the condition that $q$ be a root of unity alone says nothing about the dimensions of all representations \cite{Lusztig,DeConciniKac,DeConciniKacProcesi,Arnaudon,Pflaum}.

\subsection{Generic graded simple algebras and forms}
Let $A$ be finite-dimensional, $G$-simple and $e$-central over an algebraically closed field.  The graded identities determine its graded isomorphism class under the hypotheses of \cite{AljadeffHaile}.  The generic graded algebra of \cite{AljadeffKarasik} parametrises graded forms and their specialisations.  On each open subset over which its homogeneous components are finite locally free of fixed rank, it yields a classifying map to $\PI^G_{\idG(A),\mathbf d}$.  The orbit substack is $B\Aut_G(A)$; the PI-restricted tangent complex detects deformations preserving the complete identity ideal, while the Azumaya open substack isolates the central-simple locus.

\end{document}